\documentclass[reqno]{amsart}
\usepackage{amssymb,latexsym}

\newtheorem{thm}{Theorem}[section]
\newtheorem{lemma}{Lemma}[section]
\newtheorem{cor}{Corollary}[section]
\newtheorem{prop}{Proposition}[section]

\theoremstyle{definition}

\theoremstyle{remark}
\newtheorem{rem}{Remark}[section]
\newtheorem{exam}{Example}[section]
\makeatletter
\@addtoreset{figure}{section}
\def\thefigure{\thesection.\@arabic\c@figure}
\def\fps@figure{h, t}
\@addtoreset{table}{bsection}
\def\thetable{\thesection.\@arabic\c@table}
\def\fps@table{h, t}
\@addtoreset{equation}{section}

\makeatother

\allowdisplaybreaks

\begin{document}

\title[Diffeomorphism groups with $H^1$ metric]
{Geometry and curvature of diffeomorphism groups with
$H^1$ metric and mean hydrodynamics}
\author[S. Shkoller]{Steve Shkoller}
\address{ {CNLS,MS-B258\\Los Alamos, NM 87545}}
\address{ {CDS\\California Institute of Technology, 107-81\\Pasadena,CA 91125}}
\subjclass{Primary 58B20, 58D05; Secondary 76E99}
\date{April 15, 1998; current version July 9, 1998; To appear in J. Func. Anal.}
\email{{shkoller@cds.caltech.edu}}
\keywords{Geodesics, Hilbert diffeomorphism groups, Lagrangian stability.}

\begin{abstract}
In \cite{HMR1}, Holm, Marsden, and Ratiu derived a new model
for the mean motion of an ideal fluid in Euclidean space given by the equation
$\dot{V}(t) + \nabla_{U(t)} V(t) - \alpha^2 [\nabla U(t)]^t \cdot \triangle U(t)  = 
-\text{grad}\ p(t)$ where $\text{div} U=0$, and $V = (1- \alpha^2 \triangle)U$.
In this model, the momentum $V$ is transported by the velocity $U$,
with the effect that nonlinear interaction between modes corresponding
to length scales smaller than $\alpha$ is negligible.  We generalize this
equation to the setting of an $n$ dimensional compact Riemannian manifold.
The resulting equation is the Euler-Poincar\'{e} equation associated with
the geodesic flow of the $H^1$ right invariant metric on ${\mathcal D}^s_\mu$, 
the group of volume preserving Hilbert diffeomorphisms of class $H^s$. 
We prove that the geodesic spray is continuously differentiable
from $T{\mathcal D}_\mu^s(M)$ into $TT{\mathcal D}_\mu^s(M)$ so that a standard
Picard iteration argument proves existence and uniqueness on a finite time
interval. Our goal in this paper is to establish the foundations for Lagrangian
stability analysis following Arnold \cite{A}.  To do so, we use submanifold
geometry, and prove that the {\it weak} curvature tensor of the right invariant
$H^1$ metric on ${\mathcal D}^s_\mu$ is a bounded trilinear map in the $H^s$ 
topology, from which it follows that solutions to Jacobi's equation exist.
Using such solutions, we are able to study the infinitesimal stability behavior
of geodesics.
\end{abstract}

\maketitle

\section{Introduction}
\label{Intro}
\subsection{Background}
The Lagrangian formalism for the hydrodynamics of incompressible ideal fluids 
considers geodesic motion on ${\mathcal D}^s_\mu(M)$, the group of all volume 
preserving Hilbert diffeomorphisms of the fluid container $M$ of class $H^s$. 
Arnold \cite{A} and Ebin and Marsden \cite{EM} showed that if $\eta(t)$  is a 
smooth geodesic of the weak $L^2$ right invariant metric in 
${\mathcal D}^s_\mu(M)$,  and if $U(t)={\dot{\eta}(t)} \circ \eta(t)^{-1}$,
then the Eulerian velocity $U(t)$ is a solution of the Euler equations
\begin{equation}\label{euler}
\begin{array}{c}
\partial _t U(t) + \nabla_{U(t)}U(t) = -\text{grad}\ p(t)\\
\text{div} U(t) =0, \ \ U(0)=U_0,
\end{array}
\end{equation}
where $p(t)$ is the pressure function completely determined by $U(t)$.

The Lagrangian stability of the solutions to (\ref{euler}) is obtained
by studying the behavior of nearby geodesics.  A flow $\eta(t)$ is stable
if all geodesics in ${\mathcal D}_\mu(M)$ with sufficiently close initial
conditions at $t=0$ remain close for all $t\ge 0$.  Thus, one must study
the curvature of ${\mathcal D}^s_\mu(M)$ as this enters the linearization of
the equations of geodesic flow.  The study of the curvature of the volume
preserving diffeomorphism group with weak $L^2$ right invariant metric
was initiated by Arnold in \cite{A}.  Therein, he computed a formula for 
the sectional curvature at the identity of a group with one-side invariant 
metric in terms of the coadjoint and adjoint action, and used this formula to show
that the sectional curvature of the volume preserving diffeomorphisms of
the flat torus is negative in `many' directions.  Using this computation,
Arnold was able to demonstrate that for an idealized model of the earth's
atmosphere, deviations of fluid particles with nearby initial conditions
grow by a factor of $10^5$ in two months, making longterm dynamical weather
forecast nearly impossible.  See the book by Arnold and Khesin
\cite{AK1} (as well as \cite{AK2}) for a detailed account.
  
This work initiated a detailed study of the geometry of the volume
preserving diffeomorphism group with $L^2$ right invariant metric. Ebin
and Marsden \cite{EM} provided the differentiable structure for the
diffeomorphism groups of Sobolev class and established the functional-analytic
foundations of study (see also \cite{E}).  Lukatskii \cite{L1,L2,L3} gave
detailed explicit computations of the curvature of the  measure-preserving
diffeomorphism group on the torus.  Misio\l ek \cite{M1,M2} and
Bao, Lafontaine, and Ratiu \cite{BLR} used submanifold geometry to compute
the sectional curvature of ${\mathcal D}^s_\mu(M)$ for arbitrary manifolds $M$.
Shnirelman \cite{S1,S2} has studied the Riemannian distance on 
${\mathcal D}_\mu$ induced by the $L^2$ metric, and obtained bounds on 
the diameter of ${\mathcal D}_\mu$.  Again, see \cite{AK1} for a comprehensive
account of all of these developments.

\subsection{Motivation for the $H^1$ metric}
Our interest is in developing the geometry of the volume preserving 
diffeomorphism group with weak $H^1$ right invariant metric and studying
the properties of its curvature operator.  We are motivated by the recently
developed models of Holm, Marsden, and Ratiu \cite{HMR1}, \cite{HMR2} for
the mean hydrodynamic motion of incompressible ideal fluids in Euclidean
space.  Their basic idea was to obtain a model which averages over small 
scale fluctuations of order $\alpha$ using an additive decomposition of
a given vector field into its mean and oscillatory components.  Following
\cite{HMS}, we generalize this procedure to diffeomorphism groups of
Riemannian manifolds where mappings are `decomposed' as opposed to
vector fields.
We shall give a detailed report of this in \cite{HKMRS} for manifolds
$M$ with boundary.  Herein, we merely outline the basic construction to motivate
our study.  To do so, we shall need some notation.

Let $\alpha \mapsto \sigma^\alpha \in C^\infty([0,1],M)$.
If $U \in C^\infty(TM)$, then $U \circ\sigma \in C^\infty( TM|_
{\text{Image}(\sigma)})$.  
$U$ is said to be parallel along $\sigma$ if $\nabla_{\sigma'}U=0$,
where $\sigma'= (d/d\alpha)|_0 \sigma^\alpha$.
We set
$\alpha\mapsto P_\alpha$ to be the unique solution of
$\nabla_{\dot\sigma} P^\nabla_\alpha= 0$, $P_0=$ Id$_{T_{\sigma(0)}M}$.
$P_\alpha$ is a linear isomorphism between $T_{\sigma(0)}M$ and
$T_{\sigma(\alpha)}M$, and is called the parallel transport along $\sigma$ 
up to time $\alpha$.

We consider a geodesic curve in ${\mathcal D}_\mu(M)$ and decompose
it into its mean $\eta(t)$ and its small scale fluctuations $\zeta^\alpha(t)$ about
the mean.  The curve $\eta^\alpha(t) = \zeta^\alpha \circ \eta (t)$ describes
the motion of the fluid and is defined such that $\eta^0(t) = \eta(t)$. 
We assume that $\eta' := (d/d\alpha)|_0 \eta^\alpha$ has mean zero,
and we Taylor expand $P_\alpha^{-1}(U \circ\eta^\alpha)$ about $\alpha =0$,
where $P_\alpha$ is the parallel transport along the curve $\alpha \mapsto
\eta^\alpha(x)$.  We use the fact that $P_\alpha^{-1} \nabla_{\eta'} U
= (d/d \alpha)[ P_\alpha^{-1} U(\eta^\alpha)]$, to obtain 
$P_\alpha^{-1}U \circ \eta^\alpha = U \circ \eta + \alpha \nabla U \cdot \eta' +
O(\alpha^2)$.  Substitution of this Taylor expansion into the kinetic 
energy followed  by a computation of  its mean gives
${\frac{1}{2}}\int_M [\langle U, U \rangle + \alpha^2 \langle \nabla U,
\nabla U \rangle] \mu + O(\alpha^4)$, where $\mu$ is the volume form on $M$. 
and where, for simplicity, we  set $\overline{\eta'\otimes \eta'} = \text{Id}$.
This is not essential as the term $\langle 
\overline{\eta' \otimes \eta'} \nabla U, \nabla U\rangle$ may also be used to 
define the $H^1$ metric at the identity.

The resulting Euler-Poincar\'{e}  equation for the $H^1$ metric provides 
a new model for the mean motion of incompressible ideal fluids given by
\begin{equation}\label{CH}
\begin{array}{c}
\dot{V}(t) + \nabla_{U(t)} V(t) - \alpha^2[\nabla U(t)]^t \cdot \triangle U(t)= 
-\text{grad}\ p(t)\\ 
V = (1- \alpha^2 \triangle)U, \\
\text{div} U=0, \ U(0)=U_0.
\end{array}
\end{equation}
We call this equation the Euler-$\alpha$ equation or the averaged Euler
equation.
Unlike the Euler equation (\ref{euler}) which conserves the $L^2$ kinetic
energy $\|u\|_{L^2}$, this model conserves the $H^1$ `kinetic' energy 
$\|u\|_{H^1}$.   
Geodesic motion of the $\alpha$-$H^1$ right invariant metric on the volume 
preserving diffeomorphism group has the following effect on solutions $U$ of
(\ref{CH}):  nonlinear interaction among modes  corresponding to scales
smaller than $\alpha$ is regularized by the inversion of the elliptic
operator $(1-\alpha^2\triangle)$, so that the behavior of the solution at 
small scales is controlled by nonlinear dispersion
instead of viscous dissipation, and an $H^1$ conservation law is preserved.
Dissipation may then be added to  (\ref{CH}) to obtain a Navier-Stokes-$\alpha$
model (see \cite{FHT} for the proof of global existence of the 
Navier-Stokes-$\alpha$ model in three dimensions as well as
bounds on the dimension of the global attractor).

\subsection{Outline}
The goal of this paper is to develop the foundations for the Lagrangian
stability analysis of equation (\ref{CH}).  For our analysis, we shall set
$\alpha =1$.  Volume preserving diffeomorphism groups on Riemannian manifolds 
equipped with the $H^1$ right invariant metric have not
previously been studied, so we begin by developing the
fundamental geometric structures.  

After computing the unique Riemannian
covariant derivative of the $H^1$ right invariant metric on the diffeomorphism
group ${\mathcal D}^s(M)$, $M$ a compact Riemannian manifold, we use the Hodge 
theorem to induce the $H^1$ Riemannian covariant
derivative on ${\mathcal D}_\mu^s(M)$.  This, in turn, provides the
geodesic spray ${\mathcal S}: T{\mathcal D}_\mu^s(M) \rightarrow
TT{\mathcal D}_\mu^s(M)$ which,  just as in the case of the Euler equations,
is continuously differentiable.   A standard Picard iteration argument may
then be used to establish the existence and uniqueness of (\ref{CH}) on
a finite time interval.   In the case that the compact manifold $M$ has a boundary,
there are two very interesting subgroups of ${\mathcal D}_\mu^s(M)$ on which
the geodesic flow of the right invariant $H^1$ metric is also $C^1$.  In
\cite{HKMRS}, we define these subgroups which take into account two different
kinds of boundary conditions that may be imposed on the Euler-$\alpha$ equations.

Having this result, we proceed to study the curvature of the right invariant 
$H^1$ connection.  We follow Misio\l ek \cite{M1} and use
basic submanifold geometry, in particular the Gauss equation, to define
the curvature on the volume preserving diffeomorphism group, thought of
as a weak submanifold (and subgroup) in the weak $H^1$ topology of the
full diffeomorphism group.  We are
able to prove that this weak curvature tensor is a bounded trilinear map
in the $H^s$ topology on $M$ for $s>{\frac{n}{2}}+2$, and hence that solutions
to the Jacobi equation exist.  We note that due to the weak metric, the 
boundedness of the curvature of the $H^1$ connection cannot be immediately
infered from the regularity of the geodesic spray.

Next, we show that, just as for the Euler equations,
pressure constant flows in directions with negative sectional curvature of the
full diffeomorphism group, imply that the sectional curvature of the volume
preserving subgroup is negative, and hence that such flows are
are Lagrangian unstable, and do not possess conjugate points.

We remark, that even if $M$ is a flat manifold such as the flat
torus ${\mathbb T}^n$, the volume preserving diffeomorphism group
${\mathcal D}_\mu^s({\mathbb T}^n)$ is not flat. 
In fact, even the curvature of the right invariant $H^1$ metric on
${\mathcal D}^s({\mathbb T}^n)$ does not vanish.  Note that this is in contrast
with the curvature of the right invariant $L^2$ metric on
${\mathcal D}^s({\mathbb T}^n)$ which does vanish.

The paper is structured as follows.
In Section 2, we describe the functional analytic setting of the geometry
of the diffeomorphism group with $H^1$ metric. In Section 3, we define the 
covariant derivative of the $H^1$ metric and prove the local well-posedness of
the geodesic equations of this $H^1$ metric on the volume preserving 
diffeomorphism group.  In Section 4, we define the curvature of
the $H^1$ metric on ${\mathcal D}_\mu(M)$, prove that it is bounded in
the strong $H^s$ topology, and establish existence and uniqueness results
for the Jacobi equation.  Finally, in Section 5, we describe the Lagrangian
instability of the mean motion of incompressible ideal fluids.

\section{Functional-Analytic Setting}
\label{sec2}
\subsection{Preliminaries}
Let $(M,\langle \cdot, \cdot \rangle)$ be a compact oriented Riemannian 
$n$ dimensional manifold without
boundary and define ${\mathcal D}^s(M)$ to be the set of all bijective 
maps $\eta:M \rightarrow M$ such that $\eta$ and $\eta^{-1}$ are of
Sobolev class $H^s$.  For $s > {\frac{n}{2}} +1$, ${\mathcal D}^s(M)$
is a $C^\infty$ infinite dimensional Hilbert manifold which, about
each $\eta$, is locally diffeomorphic to the Hilbert space
$H^s_\eta(TM):= \{ X \in H^s(M,TM) : \pi \circ X = \eta\}$ where
$\pi:TM \rightarrow M$.  The condition $s > {\frac{n}{2}} +1$ ensures
that ${\mathcal D}^s(M) \subset H^s(M,M)$ is open (see \cite{MEF},
Proposition 2.3.1).

A local chart is given by $\omega_{\text{exp}}:
H^s_\eta(TM) \rightarrow {\mathcal D}^s(M)$, $\omega_{\text{exp}}(X)
= \text{exp} \circ X$, where exp is the Riemannian exponential map of
$\langle \cdot, \cdot \rangle$.  
The manifold ${\mathcal D}^s(M)$ is a topological
group with composition being the group operation.  The $\omega$-lemma
asserts that for each $\eta \in {\mathcal D}^s(M)$, right composition
$\alpha_\eta:{\mathcal D}^s(M) \rightarrow {\mathcal D}^s(M)$ is $C^\infty$,
while for all $\eta \in {\mathcal D}^{s+r}(M)$, left composition
$\omega_\eta: {\mathcal D}^s(M) \rightarrow {\mathcal D}^s(M)$ is $C^r$.

\subsection{Weak $L^2$ structure}
The weak $L^2$ right invariant Riemannian metric on ${\mathcal D}^s(M)$ is 
given by
\begin{equation}\label{1}
\langle X_\eta, Y_\eta \rangle_0 = \int_M \langle X_\eta(x), Y_\eta(x)
\rangle_{\eta(x)} \mu(x),
\end{equation}
where $\eta \in {\mathcal D}^s(M)$, $X_\eta, Y_\eta \in T_\eta{\mathcal D}^s(M)$,
and $\langle \cdot, \cdot \rangle$ and $\mu$ are the Riemannian metric and
volume element on $M$.  We let $\nabla$ be the Levi-Civita covariant derivative
of $\langle \cdot, \cdot \rangle$ on $M$, and $K:T^2M \rightarrow TM$ the
induced connector.  
\begin{rem}\label{cd}
Associated to the unique Riemannian connector $K$ of the metric 
$\langle \cdot, \cdot \rangle$ on $M$ are unique local connection $1$-forms which
which can also be used to define $\nabla$.  Let us denote by ${\mathcal V}$
the model space of $TM$. By definition, 
there exists an open cover $\{{\mathcal O}_a\}$ of $M$
and functions $ \{\psi_a\}$ defined on ${\mathcal O}_a$ such that for all
$x \in {\mathcal O}_a$, $\psi_a (s): {\mathcal V} \rightarrow T_xM$ is an
isomorphism and the map $x \mapsto \psi_a(x) \xi$ from ${\mathcal O}_a$ to
$TM$ is smooth for all $\xi \in {\mathcal V}$.
If $U \in C^\infty(TM)$ and $V \in T{\mathcal O}_a$, then
$U(x) = \psi_a(x) \xi(x)$  where  $\xi(x) = \psi_a(x)^{-1} U(x) 
\in  {\mathcal V}$ for all $x \in {\mathcal O}_a$, and
$\nabla$ on $TM$ necessarily has the form
$\nabla_V U = \psi_a(x) [ T \xi \cdot V + {\mathcal A}^a \langle V \rangle
\xi(x)]$, where the local connection $1$-forms ${\mathcal A}^a$ are defined by
${\mathcal A}^a \langle V \rangle \xi := \psi_a(x)^{-1} \nabla_V
[ \psi_a(x) \xi]$ for all $\xi \in {\mathcal V}$.
\end{rem}

It is a fact that the unique Levi-Civita $L^2$ covariant derivative $\nabla^0$
of $\langle \cdot, \cdot \rangle_0$ 
is given pointwise by $\nabla$ (see \cite{EM});  namely, if $X,Y \in
C^\infty(T{\mathcal D}^s(M))$, then
\begin{equation}\label{0}
\nabla^0_XY = K \circ (TY \cdot X).
\end{equation}
Furthermore, $\nabla^0$ is right invariant.  For $X_\eta, Y_\eta \in
T_\eta {\mathcal D}^s(M)$, let $X,Y$ be their $C^\infty$ extensions to
vector fields on ${\mathcal D}^s(M)$.  
Let $t \mapsto \eta_t$ be a smooth curve in ${\mathcal D}^s(M)$ such
that $\eta_0 = \eta$ and $(d/dt)|_0\eta_t = X_\eta$.  Then
\begin{eqnarray*}
\nabla^0_XY(\eta) &=& \left.\frac{d}{dt}\right|_0 Y(\eta_t) + \Gamma_\eta
(X_\eta,Y_\eta) \\
&=& \left.\frac{d}{dt}\right|_0 Y(\eta_t \circ \eta^{-1} )\circ \eta +
(\nabla _{X_\eta \circ \eta^{-1}} Y_\eta \circ \eta^{-1}) \circ \eta,
\end{eqnarray*}
where $\Gamma_\eta: T_\eta{\mathcal D}^s(M) \times T_\eta{\mathcal D}^s(M)
\rightarrow T_\eta {\mathcal D}^s(M)$ is the Christoffel map.  Namely,
for fixed $\eta \in {\mathcal D}^s(M)$,
let $({\mathcal O}_a, \psi_a)$ be a local frame (or trivialization) for
the bundle 
$${\mathcal E}_\eta =  \cup_{x\in M} T_{\eta(x)}M \downarrow \eta(M)$$
modeled on ${\mathcal W}$.  Then for
each $x \in {\mathcal O}_a$, $\psi_a(x): {\mathcal W} \rightarrow T_{\eta(x)}M$
is an isomorphism.  Letting $\xi(x)= \psi_a(x)^{-1}Y_\eta(x)$, for each
$x \in {\mathcal O}_a$, the Christoffel map is given by
$\Gamma_\eta(X_\eta,Y_\eta)(x)= \psi_a(x) [ {\mathcal A}^a(\eta(x)) \langle
X_\eta(x)\rangle \xi(x)]$.  The covariant derivative $\nabla$ on 
${\mathcal E}_\eta$ is given by the operator
$\nabla: C^\infty({\mathcal E}_\eta) \times {\mathcal E}_\eta \rightarrow
C^\infty({\mathcal E}_\eta)$, or for $X_{\eta(x)}$, $Y_{\eta(x)}$ elements
of the fiber ${\mathcal E}_{\eta(x)}$ over $\eta(x)$, 
$\nabla_{X_{\eta(x)}}Y_{\eta(x)} \in {\mathcal E}_{\eta(x)}$.  
It is clear that this is equivalent to
$\nabla_{(Y_\eta \circ \eta^{-1})} (X_\eta \circ \eta^{-1})\circ \eta$ 
using the symbol $\nabla$ here to denote the covariant derivative on $M$ (or $TM$).
We shall use the symbol $\nabla$ to denote the covariant derivative on  both 
$TM$ and ${\mathcal E}_\eta$, as the context will be clear.

We may also consider $M$ as the base manifold, in which case we define
the pull-back bundle $\eta^*(TM) = \cup_{x\in M} T_{\eta(x)}M \downarrow M$.
The covariant derivative on this bundle is the operator
$\nabla: C^\infty({\mathcal E}_\eta) \times TM \rightarrow
C^\infty({\mathcal E}_\eta)$.  In this setting, we differentiate a vector
$Y_{\eta(x)}$ in the direction of a vector in $TM$, and this vector is
often obtained by the push-forward of a vector $X_{\eta(x)}\in T_{\eta(x)}M$
by $\eta^{-1}$.  For example, $\nabla_{T\eta^{-1}(\eta(x)) X_\eta(x)}Y_{\eta(x)}
\in  T_\eta(x) M$.  It is often convenient for computations to take this
equivalent point of view.

\subsection{The Laplacian}
Letting $\triangle= d \delta + \delta d$ denote the Laplace-de Rham 
operator\footnote{
We identify vector fields and $1$-forms on $M$.}, we define the $H^s$
metric as follows. Let $X,Y \in T_e {\mathcal D}^s(M)$ and set
\begin{equation}\label{2}
\langle X, Y \rangle_s = \int_M \langle X(x), (1+ \triangle^s)Y(x)
\rangle \mu(x).
\end{equation}
Extending $\langle \cdot, \cdot \rangle_s$ to ${\mathcal D}^s(M)$ by right
invariance gives a smooth invariant metric on ${\mathcal D}^s(M)$.  We shall
be particularly interested in the metric $\langle \cdot, \cdot \rangle_1$.

In order to obtain formulas for the unique Levi-Civita covariant derivative  of 
$\langle \cdot, \cdot \rangle_1$, it is convenient to express the metric
(\ref{2}) in terms of the rough Laplacian $\hat{\triangle}= \text{Tr}\nabla\nabla$.
We will need the relationship between the rough Laplacian and the Laplace-de
Rham operator so that we may express (\ref{2}) in terms of $\hat{\triangle}$.
Let $\nabla^*$ denote the $L^2$ formal adjoint of $\nabla$ so that for 
any $X\in C^\infty(TM)$ and $S,T\in C^\infty(E)$, $E$ a vector bundle over
$M$, $\langle \nabla^*_X S(x), T(x)\rangle_0 = \langle S(x), \nabla_X
T(x)\rangle_0$. Then $\nabla^*_X = - \nabla_X + \text{div}X$. To see this, note
that
\begin{eqnarray*}
\langle \nabla^*_X S,T \rangle_0 = \int \langle S , \nabla_XT\rangle\mu &=&
\int X \langle S,T \rangle \mu - \langle\nabla_XS,T\rangle_0\\
&=& \int \langle S,T \rangle\text{div}X \mu- \langle\nabla_XS,T\rangle_0.
\end{eqnarray*}
If $\text{div}X=0$, then $\nabla^*_X = - \nabla_X$ which we shall often make
use of.

Next, let $\tau \in C^\infty(T^*M\otimes TM)$, let $\{e_i\}$ be a local
orthonormal frame on $M$, and let $\sigma \in C^\infty(TM)$ with support
in the domain of definition of the local frame $\{ e_i\}$.  Then
$$\langle \nabla^* \tau , \sigma \rangle_0 = 
\langle \tau , \nabla\sigma \rangle_0 = \langle \tau\langle e_i\rangle, \nabla_{e_i}
\sigma \rangle_0 = \langle \nabla_{e_i}^* ( \tau \langle e_i \rangle ), \sigma\rangle_0.$$
 We may choose the frame $\{e_i\}$, so that locally $\nabla e_i=0$
and hence $\text{div}e_i =0$.  Then
$$ \nabla^* \tau = \nabla^*_{e_i}\tau \langle e_i \rangle = -\nabla_{e_i}
( \tau \langle e_i \rangle) = - (\nabla_{e_i} \tau)\langle e_i \rangle
= -\nabla \tau (e_i,e_i),$$  
where the last equality follows from our choice of frame, since
$\nabla_{e_i} ( \tau\langle e_i \rangle) = (\nabla_{e_i}\tau) \langle e_i \rangle
= \nabla\tau \langle e_i,e_i\rangle$.  Hence $\nabla^* \tau = -\nabla \tau(e_i,e_i)$,
and since $\nabla X \in C^\infty(T^*M \otimes TM)$, we have that
$$ \hat{\triangle} = -\nabla^* \nabla.$$
With the notation established, we write Bochner's formula relating $\hat{\triangle}$
with $\triangle$ on $1$-forms as
\begin{equation}\label{5}
\triangle \alpha = \hat{\triangle} \alpha + \alpha\langle Ric \langle 
\cdot \rangle \rangle,
\end{equation}
where $Ric\langle X \rangle := R(e_i,X)e_i$, $R$ being the curvature of
$\nabla$ on $M$ (see, for example, \cite{R}).
Because the Ricci tensor is a self-adjoint operator with respect to the
metric on $TM$, for $X \in C^\infty(TM)$, we have that
$$ \triangle X = \nabla^* \nabla X + Ric\langle X \rangle.$$ 

\subsection{Weak $H^1$ metric}
Using (\ref{2}), the $H^1$ metric at the identity may be re-expressed as
\begin{eqnarray}
\langle X, Y \rangle_1 &=& \langle X,(1+ Ric)Y \rangle_{L^2}+
            \langle  X, \nabla^*\nabla Y \rangle_{L^2} \nonumber\\
        &=& \langle X,(1+ Ric)Y \rangle_{L^2}+
            \langle  \nabla X, \nabla Y \rangle_{L^2} \label{2a}
\end{eqnarray}
for all $X,Y \in T_e {\mathcal D}_\mu^s(M)$.  
The metric (\ref{2a}) extends smoothly by
right translation in the following way. Let $X_\eta,Y_\eta \in
T_\eta{\mathcal D}_\mu^s(M)$.  Then
\begin{eqnarray}
\langle X_\eta, Y_\eta \rangle_1 &=&
\int_{M} \langle X_\eta(x), Y_\eta(x) + Ric \langle Y_\eta \circ \eta^{-1}
 \rangle \circ \eta(x) \rangle_{\eta(x)}\nonumber\\
&& + \langle  \nabla(X_\eta \circ \eta^{-1}) \circ \eta(x), 
\nabla (Y_\eta \circ \eta^{-1}) \circ \eta(x) \rangle_{\eta(x)}  \mu.
\label{2b}
\end{eqnarray}

From the implicit function theorem, the set of all volume preserving $H^s$
diffeomorphisms of $M$, ${\mathcal D}_\mu^s(M):=\{ \eta \in
{\mathcal D}^s(M): \eta^*(\mu)=\mu\}$, is a submanifold of ${\mathcal D}^s(M)$
with the induced right invariant $H^1$ Riemannian metric,
as well as a subgroup.  For each $\eta \in {\mathcal D}_\mu^s(M)$, the metric
(\ref{2b}) defines a smooth orthogonal projection $P_\eta : T_\eta
{\mathcal D}^s(M) \rightarrow T_\eta {\mathcal D}_\mu^s(M)$ defined by
\begin{equation}
P_\eta(X) = (P_e (X \circ \eta^{-1})) \circ \eta, \ \ X\in 
T_\eta{\mathcal D}^s(M), \nonumber
\end{equation}
where $P_e$ is the $H^1$ orthogonal projection onto the $1$-forms $\{\alpha \in H^s:
\alpha \in \text{ker}\delta\}$ in the Hodge decomposition
\begin{equation}\label{hodge}
H^s(T^*M) = \text{ker}\delta \oplus_{H^1} dH^{s+1}(M).
\end{equation}
See \cite{Mor} for a detailed proof  of the Hodge decomposition.  
\begin{rem} We remark
here that it is essential to use the Laplace-de Rham operator in defining
the metric (\ref{2b}) in order for the Hodge decomposition to hold.  Using
the rough Laplacian instead to define the $H^1$ metric would not provide
an orthogonal decomposition in the $H^1$ topology of divergence-free vector 
fields and gradients
of functions, unless the manifold $M$ is either flat or Einstein, as can
be seen from (\ref{5}).
\end{rem}

\section{$H^1$ covariant derivative and its geodesic flow}
\label{sec3}

\subsection{Weak $H^1$ Riemannian connection}
Next, we compute the Riemannian covariant derivative on
${\mathcal D}^s(M)$ of the
$H^1$ right invariant metric restricted to vectors tangent
to ${\mathcal D}_\mu^s(M)$.  Using the Hodge decomposition,
we define the induced covariant  derivative $\tilde{\nabla}^1$ 
on ${\mathcal D}_\mu^s(M)$.
We then prove the local well-posedness of the geodesic equations
of $\tilde{\nabla}^1$.

\begin{thm}
The unique Levi-Civita covariant derivative $\nabla^1$ of 
$\langle \cdot, \cdot \rangle_1$  restricted  to vector fields
in $T{\mathcal D}_\mu^s(M)$ is given by
\begin{equation}\label{3}
\nabla^1_XY = \nabla^0_XY + A(X,Y)+B(X,Y)+C(X,Y),
\end{equation}
where  for any $\eta \in {\mathcal D}_\mu^s(M)$,
\begin{eqnarray}
&&A_\eta(X_\eta,Y_\eta)={\frac{1}{2}} (1+Ric_\eta-\hat{\triangle}_\eta)^{-1}
\left[ \nabla^* \{ \nabla X_\eta [T\eta]^{-1} \nabla Y_\eta
[T\eta]^{-1} {[T\eta]^{-1}}^t \right.\nonumber\\
&&\qquad + \nabla Y_\eta [T\eta]^{-1} \nabla X_\eta [T\eta]^{-1} {[T\eta]^{-1}}^t
+ (\nabla X_\eta [T\eta]^{-1} ) (\nabla Y_\eta [T\eta]^{-1})^t {[T\eta]^{-1}}^t\nonumber\\
&&\qquad +(\nabla Y_\eta [T\eta]^{-1}) (\nabla X_\eta [T\eta]^{-1})^t {[T\eta]^{-1}}^t
-(\nabla X_\eta [T\eta]^{-1})^t(\nabla Y_\eta [T\eta]^{-1}) {[T\eta]^{-1}}^t\nonumber\\
&&\qquad - \left.(\nabla Y_\eta [T\eta]^{-1})^t(\nabla
X_\eta [T\eta]^{-1}) {[T\eta]^{-1}}^t\} \right],\nonumber\\
\nonumber\\
&&B_\eta(X_\eta,Y_\eta)= \frac{1}{2}
\left.(1+ Ric_\eta-\hat{\triangle}_\eta)^{-1}\right\{
-\mathrm{Tr}
[ R( \nabla X_\eta T\eta^{-1}\langle \cdot \rangle , Y_\eta)\cdot  \nonumber\\
&& \qquad + R( \nabla Y_\eta T\eta^{-1}\langle \cdot \rangle , X_\eta)\cdot 
+R(X_\eta, \cdot) \nabla Y_\eta T\eta^{-1}\langle \cdot \rangle
+R(Y_\eta, \cdot) \nabla X_\eta T\eta^{-1}\langle \cdot \rangle
]
\nonumber \\
&& \qquad \left. + \nabla^*
[R( X_\eta,{T\eta^{-1}}^t)Y_\eta+
R( Y_\eta,{T\eta^{-1}}^t)X_\eta] \right\}, \nonumber \\
\nonumber\\
&&C_\eta(X_\eta,Y_\eta)= (1+ Ric_\eta-\hat{\triangle}_\eta)^{-1}
\left[ (\nabla_{X_\eta}Ric)\langle Y_\eta \rangle +
    (\nabla_{Y_\eta}Ric)\langle X_\eta \rangle \right. \nonumber \\
&&\qquad \left.
-{\frac{1}{2}}\left[
\langle (\nabla Ric\langle \cdot \rangle \langle X_\eta\rangle , Y_\eta \rangle^\sharp
+\langle (\nabla Ric\langle \cdot \rangle \langle Y_\eta\rangle , 
X_\eta \rangle^\sharp \right] -Ric_\eta\langle[X_\eta,Y_\eta]\rangle
\right],
\label{4}
\end{eqnarray}
where $X_\eta,Y_\eta \in T_\eta{\mathcal D}^s_\mu(M)$,
$$Ric_\eta \langle X_\eta \rangle =
Ric \langle X_\eta \circ \eta^{-1}\rangle \circ \eta $$
is the right-translated Ricci tensor, 
$$\hat{\triangle}_\eta=-\nabla^*[\nabla(\cdot ) (T\eta)^{-1}
{(T\eta)^{-1}}^t],$$
and $(\cdot)^\sharp$ is the operator mapping $1$-forms to vector
fields through the given metric on $M$.
\end{thm}
\begin{proof}
Formula (\ref{3}) is obtained by a lengthy computation using
(\ref{2b}) and the fundamental theorem of Riemannian geometry which
associates to every strong metric, a unique Levi-Civita covariant 
derivative.  Although $\langle \cdot, \cdot \rangle_1$ is a weak
metric, $\nabla^1$ is still uniquely defined by virtue of the existence of
a $C^1$ geodesic spray restricted to tangent vectors on ${\mathcal D}_\mu^s(M)$
(see Theorem \ref{thm_spray}).
\end{proof}

\begin{rem}
Note that for $X_\eta \in H^s_\eta(TM)$,  the operators $[T\eta]^{-1}$, 
${[T\eta]^{-1}}^t$, and $\nabla X_\eta$ induce the  following pointwise
operators
\begin{equation}\nonumber
\begin{array}{c}
[T\eta(x)]^{-1} : T_{\eta(x)}M \rightarrow T_xM, \\
 {[T\eta(x)]^{-1}}^t  :  T_xM \rightarrow T_{\eta(x)}M, \\
(\nabla X_\eta)(x) : T_xM \rightarrow T_{\eta(x)}M.
\end{array}
\end{equation}
\end{rem}

\begin{rem}
Since $[T\eta]^{-1}{[T\eta]^{-1}}^t$ is positive symmetric, the spectrum of
$-\hat{\triangle}_\eta$, $\sigma(-\hat{\triangle}_\eta)$, is positive.
We can ensure that $0\not\in \sigma(1+Ric_\eta-\hat{\triangle})$ by
requiring that $M$ have nonnegative Ricci curvature or in the case that $M$ has
negative Ricci curvature, by insisting that $|-\sigma(Ric_\eta)| \le 1$.
More generally, we require Ker$(1 + Ric_\eta - \hat{\triangle}_\eta)$ to 
be either empty or unique for all $x \in M$, $\eta \in {\mathcal D}_\mu^s(M)$.
In the case that the kernel is not empty, 
we shall restrict our phase space to the orthogonal
complement of Ker$(1 + Ric_\eta - \hat{\triangle}_\eta)$ but this may only occur
if on manifolds $M$ with negative Ricci curvature (this is essentially Bochner's
theorem).
\end{rem}

Now, on $H^{s+1}(M)$, $\triangle= d \delta = -\text{div grad}$, so an explicit
formula for $P_e:T_e {\mathcal D}^s(M) \rightarrow T_e {\mathcal D}_\mu^s(M)$
is obtained  as follows.  Suppose that $V \in H^s(TM)$, and
let $p\in H^{s+1}(M)$ solve $\triangle p = \text{div}V$.  Then
$$P_e(V) = V - \text{grad} \triangle^{-1} \text{div}V.$$

We shall denote the orthogonal projection onto $dH^{s+1}(M)$ by
\begin{equation}\label{Qe}
Q_e(V)=\text{grad} \triangle^{-1} \text{div}V.
\end{equation}

${\mathcal D}_\mu^s(M)$ thus becomes a weak Riemannian submanifold of
${\mathcal D}^s(M)$ with the metric (\ref{2b}), and the induced covariant
derivative 
$$\tilde{\nabla}^1 = P \circ \nabla^1$$ 
is inherited from ${\mathcal D}^s(M)$.  
\subsection{Geodesic flow of $\tilde{\nabla}^1$}
\begin{thm}\label{thm_ep}
If $\eta(t)$ is a geodesic of $\tilde{\nabla}^1$, 
then $U(t) = \dot{\eta} \circ \eta^{-1}(t)$ is a vector field on $M$ which 
satisfies the mean motion equations of an ideal fluid,
\begin{equation}\label{7}
\begin{array}{c}
\partial_t U(t) + (1+{\triangle})^{-1}\left[\nabla_{U(t)} 
(1+{\triangle})U(t) + \langle \nabla U(t)\langle \cdot  \rangle,
{\triangle} U(t) \rangle^\sharp \right] = -\mathrm{grad}\ p(t) \\
\mathrm{div}U(t)=0, \ U(0)= U_0, 
\end{array}
\end{equation}
where $p(t)$ is the pressure function which is determined from $V(t)$.
Laplacian 
\end{thm}
\begin{proof}
Together with the Hodge decomposition (\ref{hodge}), a
a straightforward computation of the coadjoint action $\text{ad}^*$ 
of ${\mathcal D}_\mu^s(M)$ given by
\begin{equation}
\begin{array}{c}
\langle \text{ad}^*_V W, U \rangle_1 = 
\langle \text{ad}_V U,W\rangle_1, \\
\text{ad}_UV = -[U,V], \ \ \ U,V,W \in T_e{\mathcal D}_\mu^s(M)
\end{array}
\end{equation}
shows that (\ref{7}) is  simply
$$\dot{U}(t) = -P_e \circ \text{ad}^*_{U(t)}U(t),$$
the Euler-Poincar\'{e} equation for the induced $H^1$ metric on 
${\mathcal D}_\mu^s(M)$.
\end{proof}
\begin{rem}
Notice that the Euler-Poincar\'{e} equation (\ref{7}) is expressed in
terms of the Laplace-de Rham operator $\triangle$.  In terms
of the rough Laplacian $\hat{\triangle}$, 
$$ \text{ad}^*_UU = P_e \circ
(1+Ric -\hat{\triangle})^{-1}\left[\nabla_U 
(1+Ric -\hat{\triangle})U - \nabla U^t \cdot [Ric
+ \hat{\triangle}] U\right]. $$
\end{rem}

We shall need the following lemmas, the first of which  is similar to Lemma 2 of Appendix
A in \cite{EM}.
 
\begin{lemma} \label{l1}
Let ${\hat\triangle}_{(\cdot)} :  \cup_{\eta \in
{\mathcal D}^s_\mu(M)} H^s_\eta(TM) \downarrow {\mathcal D}^s_\mu(M)
\longrightarrow \cup_{\eta \in
{\mathcal D}^s_\mu(M)} H^{s-2}_\eta(TM) \downarrow {\mathcal D}^s_\mu(M)$
be given by
$${\hat\triangle}_\eta= -\nabla^*[\nabla( \cdot) (T\eta)^{-1}(T\eta)^{{-1}^t}]$$
and the identity on ${\mathcal D}^s_\mu(M)$.  Then
${\hat\triangle}_{(\cdot)}$ is a $C^1$ bundle map.
\end{lemma}
\begin{proof}
Let $H^{s-1}_\eta(T^*M\otimes TM) = H^{s-1}(\cup_{x\in M} (T^*_{\eta(x)}M
\otimes T_{\eta(x)}M) \downarrow M)$, and
let $$f(\eta)= \nabla(\cdot) (T\eta)^{-1} (T\eta)^{{-1}^t}.$$
 
We first show that $f$ is a $C^1$ section of the bundle
$$\cup_{\eta \in {\mathcal D}^s_\mu(M)}
{\mathrm Hom}(H^s_\eta(TM), H^{s-1}_\eta(T^*M \otimes TM)) \downarrow
{\mathcal D}^s_\mu(M).$$
 
Continuity of $f$ is clear.  We compute its
derivative.  With $V \in H^s_\eta(TM)$, the $\omega$-lemma asserts that
$$Df(\eta) \langle V \rangle = \nabla( \cdot ) [T\eta]^{-1}(\nabla V)
[T\eta]^{-1}{[T\eta]^{-1}}^t - \nabla(\cdot)[T\eta]^{-1}{[T\eta]^{-1}}^t
(\nabla V)^t {[T\eta]^{-1}}^t.$$
 
Now,
\begin{eqnarray*}
&&\| Df(\eta)\|_{{\mathcal L}(H^s_\eta(TM),\mathrm{Hom}
( H^s_\eta(TM), H^{s-1}_\eta(T^*M \otimes TM))   )} \\
&& \qquad = \sup_{V\in H^s_\eta(TM), \|V\|_s=1}
\|Df(\eta)\langle V \rangle \|
_{\mathrm{Hom}(H^s_\eta(TM), H^{s-1}_\eta(T^*M \otimes TM))} \\
&& \qquad = \sup_{V\in H^s_\eta(TM), \|V\|_s=1}\sup_{W\in H^s_\eta(TM), \|W\|_s=1}
\|(Df(\eta)\langle V \rangle)\langle W \rangle \|
_{H^{s-1}_\eta(T^*M \otimes TM)} \\
&& \qquad \le C(\|T\eta\|_{s-1}, \|[T\eta]^{-1}\|_{s-1}< \infty,
\end{eqnarray*}
where the last two inequalities are due to the $\omega$-lemma and the
fact that $[T\eta]^{-1} \in H^{s-1}$ whenever $\eta \in H^s$, again by
the $\omega$-lemma.  Let ${\mathcal O} \subset {\mathcal D}^s_\mu(M)$ be a 
be neighborhood of some $\eta$.  Locally ${\hat\triangle}_{\cdot}$ acts
on ${\mathcal O} \otimes {\mathcal F}$ for a trivialization $\{\psi(\eta)\}_{\eta 
\in {\mathcal O}}$ such that $\psi(\eta): H^s_\eta(TM) \rightarrow {\mathcal F}$
isomorphically.

Computing the supremum of
$$ \| Df(\eta)\|_{{\mathcal L}(H^s_\eta(TM),\mathrm{Hom}
( H^s_\eta(TM), H^{s-1}_\eta(T^*M \otimes TM))   )}$$
over all $\eta \in {\mathcal O}$  defines the $C^1$ topology. 
Since we may bound the supremum, we have 
proven that $f$ is $C^1$.  Now
thinking of $\nabla (\cdot) [T\eta]^{-1} {[T\eta]^{-1}}^t$ as a map on
${\mathcal F}$, it  is smooth by the $\omega$-lemma.  
To see this, it suffices to consider the fiber over the identity $e$,
where the operator is a linear and hence a smooth bundle map.

The operator $\nabla^*$ acts fiberwise, and is linear, hence smooth as a bundle
map.  This proves that ${\hat\triangle}_{(\cdot)}$ is a $C^1$ bundle map,  
which proves the lemma.
\end{proof}
\begin{rem}
Although we shall only need the $C^1$ regularity, it seem likely that
by considering higher order derivatives of $\nabla (\cdot) [T\eta]^{-1}
{[T\eta]^{-1}}^t$, thought of as a bundle map, we could obtain the 
$C^k$ regularity of $\hat{\triangle}_{( \cdot )}$ for any nonnegative integer 
$k$.
\end{rem}

\begin{lemma} \label{l2}
The operator $(1 +Ric_{(\cdot)}- {\hat\triangle}_{(\cdot)})^{-1} : 
\cup_{\eta \in {\mathcal D}^s_\mu(M)} H^s_\eta(TM) \downarrow {\mathcal D}^s_\mu(M)
\longrightarrow
\cup_{\eta \in {\mathcal D}^s_\mu(M)} H^s_\eta(TM) \downarrow {\mathcal D}^s_\mu(M)$
is a $C^1$ bundle map.
\end{lemma}
\begin{proof}
By the smoothness of right translation, the map $\eta \mapsto Ric_\eta$ is smooth.
Thus, $(1 +Ric_{(\cdot)}- {\hat\triangle}_{(\cdot)})$ is smooth (using Lemma \ref{l1})
and by assumption has trivial kernel and closed range, hence is a $C^1$ bijection.
By the inverse function theorem, a $C^1$ bijective bundle map covering the identity
has a $C^1$ inverse.
\end{proof}

For the following theorem, recall that $TT{\mathcal D}_\mu^s(M)$ is
identified with $H^s$ maps ${\mathcal Y}:M \rightarrow TTM$ covering
some $X_\eta \in T_\eta{\mathcal D}_\mu^s(M)$.

\begin{thm}\label{thm_spray}
For $s>{\frac{n}{2}}+1$, there exists a neighborhood of $e \in 
{\mathcal D}_\mu^s(M)$ and an $\epsilon>0$ such that for any $V
\in T_e{\mathcal D}_\mu^s(M)$ with $\|V\|_s < \epsilon$, there exists
a unique geodesic $\dot\eta \in C^1 ((-2,2) ,T{\mathcal D}_\mu^s(M))$ satisfying
$$ \tilde{\nabla}^1_{\dot \eta}\dot\eta =0, \ \ \
\eta(0)=e,\  \dot{\eta}(0)=V,$$
with smooth dependence on $V$.
\end{thm}
\begin{proof}
Let $\eta(t)$ be a curve in ${\mathcal D}_\mu^s(M)$.
Using the formula for the induced covariant derivative of the $H^1$ metric (\ref{3})
on ${\mathcal D}_\mu^s(M)$ or by a computation of the first variation of
the energy (see \cite{HKMRS} for the detailed computation)
\begin{equation}\label{energy}
{\mathcal E}(\eta)={\frac{1}{2}}\int_{\mathbb R} \langle \dot\eta(t), \dot\eta(t) \rangle_1dt,
\end{equation}
we find that
\begin{eqnarray}
&&P_\eta \circ \nabla_{\dot{\eta}}{\dot{\eta}}= 
P_\eta \circ (1+Ric_\eta - \hat{\triangle}_\eta)^{-1}
%\left[
\bigl[
\nabla^*\left[\left\{ -(\nabla{\dot{\eta}} [T\eta]^{-1})^t (\nabla {\dot{\eta}}
%[T\eta]^{-1}) \right.\right.\right.\nonumber\\
[T\eta]^{-1}) \right.\right.\nonumber\\
&&\qquad \qquad \left. + \nabla{\dot{\eta}} [T\eta]^{-1}
\nabla{\dot{\eta}} [T\eta]^{-1} 
 + \left. (\nabla{\dot{\eta}}[T\eta]^{-1}) ( \nabla{\dot{\eta}}[T\eta]^{-1})^t
\right\} {[T\eta]^{-1}}^t \right]\nonumber \\
&&\qquad \qquad + (\nabla_{\dot{\eta}} Ric) \langle {\dot{\eta}} \rangle -
{\frac{1}{2}}\langle \nabla Ric\langle \cdot \rangle \langle \dot{\eta} \rangle,
 \dot{\eta} \rangle^\sharp
- \{\text{Tr}[ R(\nabla \dot{\eta} T\eta^{-1} \langle \cdot \rangle,
\dot{\eta}) \cdot \nonumber\\
&&\left. \qquad \qquad + R(\dot{\eta}, \cdot)\nabla \dot{\eta} T\eta^{-1} 
\langle \cdot \rangle]
+ \nabla^*\{ R(\dot{\eta}, {T\eta^{-1}}^t){\dot{\eta}} \} \right]  \label{spray}\\
&&\qquad \qquad := P_\eta \circ F_\eta({\dot{\eta}}) \nonumber.
\end{eqnarray}

Using the notation of Remark \ref{cd}, we let $({\mathcal O}_a, \psi_a)$ be a
trivialization of ${\mathcal E}_\eta$ and set ${\dot{\eta}}(x) = \psi_a(x)\xi(x)$.
For all $x \in {\mathcal O}_a$,
we express $\nabla_{\dot{\eta}(x)}{\dot{\eta}(x)}$ by
$\nabla_{\dot{\eta}}{\dot{\eta}}(x)=\psi_a(x) 
[\dot{\xi}+({\mathcal A}^a \circ \eta )(x) \langle {\dot{\eta}} \rangle
\xi(x)]$.   Let $\tilde F_\eta$ be the localization of $F_\eta$ in 
$({\mathcal O}_a, \psi_a)$.  Then, in this trivialization, we may write
(\ref{spray}) in the form of a geodesic spray ${\mathcal S}:T{\mathcal D}_\mu^s(M)
\rightarrow TT{\mathcal D}_\mu^s(M)$.  We have, locally, that
$$ {\mathcal S}_\eta({\dot{\eta}}) = \frac{d}{dt}(\eta , \psi_a^{-1} {\dot{\eta}})=
(\xi, Q_\eta \psi_a \dot \xi - P_\eta[\psi_a({\mathcal A}^a\circ \eta)   \langle
 \psi_a \xi\rangle \xi - \psi_a \tilde F_\eta]).$$
We show that ${\mathcal S}_\eta$ is a quadratic form. Clearly, $F_\eta$ is quadratic;
as for the term $Q_\eta \psi_a \dot\xi$, we note that 
$$ \dot \xi = \psi_a^{-1} \left[ \left(\psi_a \xi\circ\eta^{-1} + \nabla_{\psi_a\xi\circ
\eta^{-1}}(\psi_a \xi\circ\eta^{-1})\right) \circ \eta \right],$$
and since div$(\psi_a \xi \circ\eta^{-1})=0$, 
$Q_e( \psi_a \dot\xi\circ \eta^{-1}) \circ \eta= Q_e[ 
T(\psi_a\xi\circ\eta^{-1}) \cdot(\psi_a\xi\circ\eta^{-1})] \circ \eta$,
so that
\begin{eqnarray*}
&&Q_e( \psi_a \dot\xi\circ \eta^{-1}) \circ \eta +
Q_\eta [\psi_a({\mathcal A}^a \circ \eta) \langle  \psi_a \xi\rangle \xi] 
=Q_e[\nabla_{\psi_a\xi\circ\eta^{-1}}(\psi_a\xi\circ\eta^{-1}] \circ \eta\\
&& = \mathrm{grad} \triangle^{-1} \left[
\mathrm{Ric}(\psi_a \xi \circ \eta^{-1},\psi_a \xi \circ \eta^{-1})
+\mathrm{Tr}(\nabla(\psi_a \xi \circ \eta^{-1})\cdot \nabla(\psi_a \xi \circ \eta^{-1}))
\right] \circ \eta,
\end{eqnarray*}
where Ric$(V,W)=Ric\langle V\rangle W$.  This shows that $S_\eta$
is quadratic in $\xi$.

The projection $P_\eta$ is a smooth bundle map. 
Namely, $P:T{\mathcal D}^s(M)\downarrow {\mathcal D}_\mu^s(M) \rightarrow
T{\mathcal D}_\mu^s(M)$ is $C^\infty$.
(To prove this one need only replace
the $L^2$ orthogonal projection onto the harmonic forms by the $H^1$
orthogonal projection onto harmonic forms in Lemma 4 of Appendix A in
\cite{EM}.)

The map $x \mapsto ({\mathcal A}^a\circ \eta)(x) \in 
C^\infty({\mathcal O}_a, [T^*_{\eta(x)}M]^2 \otimes T_{\eta(x)}M)$ since
the local connection $1$-forms and right translation are both smooth maps.
Since $\psi_a(x)$ is an isomorphism, $\psi_a[({\mathcal A}^a 
\circ \eta )\langle \cdot \rangle (\cdot)]: (H^s_\eta)^2 \rightarrow 
H^s_\eta$ smoothly.

By Lemma \ref{l2},
$(1+Ric_{(\cdot)}-\hat{\triangle}_{(\cdot)})^{-1}$ is a $C^1$ bundle
map.
Since $R$ and $Ric$ are fiberwise multilinear maps, it follows from
the smoothness of right translation that 
all terms involving the curvature are smooth bundle maps.
Letting $U={\dot\eta}\circ\eta^{-1}$, we need only prove that the
terms $[-(\nabla U)^t (\nabla U) + (\nabla U)(\nabla U)+(\nabla U)
(\nabla U)^t] {[T\eta]^{-1}}^t$ are $C^1$ bundle maps.  The argument
for this is identical to that of Lemma \ref{l1}.

We have shown that $S:T{\mathcal D}^s_\mu(M) \rightarrow TT {\mathcal D}^s_\mu(M)$
is a $C^1$ bundle map.
A standard Picard iteration argument for ordinary differential equations in a 
Banach space then proves the existence of a unique $C^1$ flow (see
\cite{La}, Theorem 1.11), and this proves the theorem.
\end{proof}
Together with Theorem \ref{thm_ep}, we have proven the local well-posedness
of the Cauchy problem for the hydrodynamic mean motion equations (\ref{7})
on $M$.  This implies the following facts. \footnote{We would like to thank
the referee for pointing these out and suggesting their inclusion in this paper.}

\begin{cor}
Let $\eta \in {\mathcal D}_\mu^s(M)$ be in a sufficiently small neighborhood
of $e$.  Then, there exists a vector field $V$ on $M$ such that exp$_e(V) =
\eta$.  In other words, the Euler-$\alpha$ flow with initial condition $V$
reaches $\eta$ in time $1$.
\end{cor}

As another corollary, we immediately have the $H^1$ analog of Theorem 12.1 of
\cite{EM}.  
\begin{cor}
For $s> \frac{n}{2}+1$, let $\eta(t)$ be a geodesic of the right invariant $H^1$ metric
on ${\mathcal D}_\mu^s(M)$. If $\eta(0) \in {\mathcal D}_\mu^{s+k}(M)$  and
${\dot{\eta}}(0) \in T_{\eta(0)}{\mathcal D}_\mu^{s+k}(M)$ 
for $0 \le k \le \infty$, then $\eta(t)$ is $H^{s+k}$
on $M$ for all $t$ for which $\eta(t)$ was defined in ${\mathcal D}_\mu^s(M)$.
\end{cor}
\noindent
The proof of this theorem exactly follows the proof of Theorem 12.1 of \cite{EM}
once we have the regularity properties of the exponential map.  As noted in \cite{EM} for
the case of the Euler equations, this has the important consequence that the time of
existence of a geodesic does not depend on $s$, so that a geodesic with $C^\infty$ 
initial conditions is a curve in 
$${\mathcal D}_\mu(M) = \cap_{s> n/2} {\mathcal D}_\mu^s(M),$$
where ${\mathcal D}_\mu(M)$ is the ILH (inverse limit Hilbert) Lie group of $C^\infty$
diffeomorphisms.

\begin{rem}
A computation of the first variation of (\ref{energy}) on the full diffeomorphism
group shows that the geodesic spray has no derivative loss in this case as well.  
For example, on ${\mathbb S}^1$, with $\triangle:= \eta_x^{-1}(\partial_x \eta_x^{-1}\partial_x)$
and for $\alpha > 0$, the principle part of the geodesic spray, for $s >5/2$, is given by
\begin{equation}\label{1D}
\ddot{\eta} = (1-\alpha^2\triangle)^{-1}\left[ 
(-2\dot\eta + \alpha^2\triangle \dot\eta)
\eta_x^{-1} {\dot\eta}_x \right].
\end{equation}
It is clear that the  nonlinear dispersion arising from the $H^1$ metric
regularizes the shock formation of the Burger-Riemann equation into traveling
peaked solitons (see \cite{HMR1}).  The fact that the Burger-Riemann equation which
arises from the $L^2$ right invariant metric shocks, is a connected to the loss of
smoothness of the spray, for in the $\alpha=0$ limit, (\ref{1D}) is 
$\ddot \eta = -2 \eta_x^{-1}{\dot{\eta}}_x {\dot{\eta}}$ which has derivative
loss.

A similar but lengthier computation shows that for $s> n/2+2$, the geodesic spray 
has no derivative loss on the full diffeomorphism group in $n$ dimensions, so that
the covariant derivative $\nabla^1$ can be uniquely defined for all vectors in
$T{\mathcal D}^s(M)$.
\end{rem}

\section{Curvature of the $H^1$ metric}
\label{sec4}
Because the Lie-theoretic computation of the sectional curvature is
difficult to compute on manifolds $M$ with nonvanishing curvature, we
use basic submanifold geometry to estimate the curvature of the 
$H^1$ metric on ${\mathcal D}_\mu(M)$ for arbitrary smooth manifolds. 

\subsection{Curvature of $\nabla^1$}
We denote by $R^0$ the curvature of the $L^2$ metric $\nabla^0$.
Proposition 3.4 of \cite{M1} states that $R^0$ is completely determined
by $R$, the curvature of $M$, and is a bounded trilinear map in the
$H^s$ topology.  Namely, for $X_\eta, Y_\eta,Z_\eta \in 
T_\eta{\mathcal D}^s(M)$ and using the right invariance of
$\nabla^0$, it is evident from formula (\ref{0}) that $R^0$ may
be expressed as 
$$R^0(X_\eta,Y_\eta)Z_\eta = (R(X_\eta \circ \eta^{-1}, Y_\eta \circ
\eta^{-1})Z_\eta \circ \eta^{-1}) \circ \eta.$$
It follows that $R^0$ is right invariant, and that
$$ \|R^0_\eta(X_\eta, Y_\eta)Z_\eta\|_s \le C
\|X_\eta\|_s \|Y_\eta\|_s \|Z_\eta\|_s,$$
where $C$ denotes any constant which may depend on $s, \eta$, and the
derivatives of the metric $\langle  \cdot, \cdot \rangle$ on $M$.

Now for each $\eta \in {\mathcal D}_\mu^s(M)$,
the weak metric (\ref{2b}) splits $T_\eta{\mathcal D}^s(M)$ into the
direct sum
$$ T_\eta{\mathcal D}^s(M) = T_\eta{\mathcal D}_\mu^s(M) \oplus_{H^1}
\nu_\eta{\mathcal D}_\mu^s(M),$$
where $\nu_\eta{\mathcal D}_\mu^s(M)$ is the $H^1$ orthogonal complement
of $T_\eta{\mathcal D}_\mu^s(M)$ in $T_\eta{\mathcal D}^s(M)$.
We now introduce the (weak) second fundamental form $S$ of 
${\mathcal D}_\mu^s(M)$ by assigning to each $\eta \in {\mathcal D}_\mu^s(M)$
a map
$$S_\eta: T_\eta{\mathcal D}_\mu^s(M) \times T_\eta{\mathcal D}_\mu^s(M)
\rightarrow \nu_\eta {\mathcal D}_\mu^s(M).$$
Given $X_\eta,Y_\eta \in T_\eta{\mathcal D}_\mu^s(M)$, we extend them to
$C^\infty$ vector fields $X,Y$ on ${\mathcal D}_\mu^s(M)$, and define
\begin{eqnarray}
S_\eta(X_\eta,Y_\eta) &=& Q_\eta(\nabla^1_XY(\eta)), \label{8}\\
 &=& Q_\eta(\nabla^0_XY(\eta) + A_\eta(X_\eta,Y_\eta)) + B_\eta(X_\eta,
 Y_\eta)+ C_\eta(X_\eta,Y_\eta)), \nonumber
\end{eqnarray}
where $\eta \in {\mathcal D}_\mu^s(M)$ and
$$ Q_\eta(X_\eta) =(Q_e(X_\eta \circ \eta^{-1})) \circ \eta$$
can be computed explicitly from (\ref{Qe}).

We next define the (weak) Riemannian curvature tensor $R^1$ of
$\langle \cdot, \cdot\rangle_1$ on ${\mathcal D}^s(M)$.  This is the
trilinear map
\begin{equation}\nonumber
\begin{array}{c}
R^1_\eta: T_\eta{\mathcal D}^s(M) \times T_\eta{\mathcal D}^s(M)
\times T_\eta{\mathcal D}^s(M) \rightarrow T_\eta{\mathcal D}^s(M),\\
R^1_\eta(X_\eta,Y_\eta)Z_\eta = (\nabla^1_X \nabla^1_Y Z)_\eta -
(\nabla^1_Y \nabla^1_X Z)_\eta - (\nabla^1_{[X,Y]}Z)_\eta,
\end{array}
\end{equation}
where $\eta \in {\mathcal D}^s(M)$ and $X,Y,Z$ are smooth extensions of
vectors $X_\eta,Y_\eta,Z_\eta$ to a neighborhood of $\eta$.

\begin{lemma} \label{lem1}
For $\eta \in {\mathcal D}^s(M)$, $B_\eta : (H^s_\eta(TM))^2 \rightarrow
H^{s+1}_\eta(TM)$ continuously.
\end{lemma}
\begin{proof}
Let $X,Y,Z \in T_e{\mathcal D}^s(M)$.  Since $s>{\frac{n}{2}}+1$, $H^r$
is a multiplicative algebra for $r\ge s-1$; hence, it suffices to
obtain the estimate at the identity $e$. 

We use the fact that
$R^0$ is a continuous trilinear map in the $H^s$ topology, and
estimate $B_\eta$ using equation (\ref{4}). For the terms 
Tr$[R(\nabla_\cdot X,Y)\cdot \ + R(\nabla_\cdot Y,X)\cdot \ +
R(X,\cdot)\nabla_\cdot Y + R(Y,\cdot)\nabla_\cdot X]$
we use the
continuous embedding $H^{s-1}(TM)\hookrightarrow C^0(TM)$, while for
the term $\nabla^*[R(X,\cdot )Y + R(Y,\cdot )X]$
we use that $\nabla^*:H^s \rightarrow
H^{s-1}$ is continuous. Since $(1-\hat{\triangle})^{-1}$ is a
pseudodifferential  operator of order -2, we obtain that
$$\|B(X,Y)\|_{s+1} \le C \|X\|_s \|Y\|_s,$$
where the constant $C$ may depend on $R$ and $s$.
\end{proof}
The same argument shows that
\begin{cor}\label{cor1}
For each $\eta \in {\mathcal D}^s(M)$, $B_\eta:
H^s_\eta(TM) \times H^{s-1}_\eta(TM)  \rightarrow H^s_\eta(TM)$ continuously.
\end{cor}
\noindent
Similarly,
\begin{lemma} \label{lem2}
For each $\eta \in {\mathcal D}^s(M)$, the following are
bounded multilinear maps:
\begin{itemize}
\item[i)]$C_\eta:(H^s_\eta(TM))^2 \rightarrow H^{s+1}_\eta(TM)$,
\item[ii)]for each $X_\eta \in T_\eta{\mathcal D}^s(M)$,
$\nabla^0_{X_\eta} : H^s_\eta(TM) \rightarrow H^{s-1}_\eta(TM)$,
\item[iii)]$A_\eta:(H^s_\eta(TM))^2 \rightarrow H^{s}_\eta(TM)$.
\end{itemize}
\end{lemma}
\begin{proof}
Items i) and ii) are trivial, while for item iii), we use that $H^{s-1}$
is a Schauder ring.
\end{proof}

\begin{prop} \label{thmR1}
Let $M$ be a compact $n$ dimensional manifold.  For $s>{\frac{n}{2}}+2$,
and $\eta \in {\mathcal D}_\mu^s(M)$,
$R^1_\eta:(T_\eta{\mathcal D}_\mu^s(M))^3 \rightarrow T_\eta
{\mathcal D}_\mu^s(M)$ is continuous in the $H^s$ topology.
\end{prop}

\begin{proof}
For $\eta \in {\mathcal D}_\mu^s(M)$,
let $X_\eta,Y_\eta,Z_\eta \in T_\eta{\mathcal D}_\mu^s(M)$,  and
let $X,Y,Z$ be smooth extensions to a neighborhood of $\eta$.
Let $D(X,Y)=A(X,Y)+B(X,Y)+C(X,Y)$.  Then
\begin{eqnarray*}
R^1_\eta(X_\eta,Y_\eta)Z_\eta &=& (\nabla^1_X \nabla^1_Y Z)(\eta)-
 (\nabla^1_Y \nabla^1_X Z)(\eta) - (\nabla^1_{[X,Y]} Z)(\eta)\\
 &=& R^0_\eta(X_\eta,Y_\eta)Z_\eta +
 D(X, \nabla^1_YZ)(\eta) - D(Y, \nabla^1_XZ)(\eta) \\
 && + (\nabla^0_X D(Y,Z))(\eta) - (\nabla^0_Y D(X,Z))(\eta)\\
 && + D(X,D(Y,Z))(\eta)- D(Y,D(X,Z))(\eta)- D([X,Y],Z)(\eta).
\end{eqnarray*}

Since $R^0$ is a bounded trilinear map in the $H^s$ topology, we must
show that the remaining terms  are bounded trilinear maps in $H^s$ as
well.  These terms are of two types. Type $I$  terms involve commutation
between $\nabla^0$ and $D$, while the type $II$ terms involve commutation
between the bilinear forms $A,B$, and $C$.  From Lemmas \ref{lem1} and
\ref{lem2} it is clear that the trilinear map formed by type $II$ terms
are bounded maps in the $H^s$ topology; hence, we estimate type $I$
terms.

We begin with type $I$ terms which are the commutation of $\nabla^0$
and $B$.  Since for each $\eta\in {\mathcal D}_\mu^s(M)$, $H^{s-2}_\eta$
is a Schauder ring, using the right invariance of $\| \cdot \|_s$
it suffices to obtain the continuity of the
trilinear maps at the identity $e$.  Using Lemma \ref{lem1},
it is clear that terms of the type $\nabla^0_{X}B(Y,Z)$ are 
continuous in $H^s$, while Corollary \ref{cor1} gives the bound on
the remaining terms involving $B$.  Clearly, since $C_\eta$ is
as regularizing as $B_\eta$, by the same argument, we have that
all type $I$ terms involving the commutation of $\nabla^0$ and $C$
are continuous trilinear maps in $H^s$ as well.  The difficult type
$I$ terms to estimate are those involving the commutation of 
$\nabla^0$ and $A$, since by part iii) of Lemma \ref{lem2}, it 
appears as though a derivative loss may occur in some of these terms.

In fact, such a derivative loss does not occur, and for 
the purpose of estimating these terms, it will suffice to replace
$A_e$ with 
$$ \bar{A}(X,Y) = \hat{\triangle}^{-1} \nabla^*(\nabla X \cdot \nabla Y)$$
for $X,Y \in T_e{\mathcal D}_\mu^s(M)$.  The terms we must estimate are
given by
\begin{gather}
\nabla_Y \hat{\triangle}^{-1} \nabla^*(\nabla X \cdot \nabla Z)
+ \hat{\triangle}^{-1}\nabla^*(\nabla Y \cdot \nabla_XZ) 
+ \hat{\triangle}^{-1} \nabla^*
(\nabla Y \cdot \hat{\triangle}^{-1}\nabla^*(\nabla X \cdot \nabla Z))
\nonumber\\
-\nabla_X \hat{\triangle}^{-1} \nabla^*(\nabla Y \cdot \nabla Z)
- \hat{\triangle}^{-1}\nabla^*(\nabla X \cdot \nabla_YZ) 
- \hat{\triangle}^{-1} \nabla^*
(\nabla X \cdot \hat{\triangle}^{-1}\nabla^*(\nabla Y \cdot \nabla Z))
\nonumber\\
-\hat{\triangle}^{-1}\nabla^*( \nabla[X,Y] \cdot \nabla Z). \label{A}
\end{gather}
We shall need the following lemma which is Corollary 4.2 of \cite{T}.
\begin{lemma}\label{pd}
Let $\alpha$ and $\beta$ be pseudodifferential operators with symbols
of order $m$ and $n$, respectively.  Then the commutator $[\alpha,
\beta]$ is a pseudodifferential operator with symbol of order
$m+n-1$.
\end{lemma}
Using Lemma \ref{pd}, $[\hat{\triangle}^{-1} \nabla^*, \nabla_Y]$ is a
pseudodifferential operator of order $-1$, so that
$[\hat{\triangle}^{-1} \nabla^*, \nabla_Y]:H^s \rightarrow H^{s+1}$
continuously.  Hence, using the property of the Schauder ring, it is
clear that
$$\|[\hat{\triangle}^{-1} \nabla^*, \nabla_Y](\nabla X \cdot \nabla Z)\|_s
\le C \|X\|_s \|Y\|_s \|Z\|_s,$$
where, in general, the constant $C$ may depend on $M$ and $\eta$.  Similarly,
we have the identical estimate for
$[\hat{\triangle}^{-1} \nabla^*, \nabla_X](\nabla Y \cdot \nabla Z)$.

Next, we consider the endomorphism
$$\nabla_Y \nabla X \cdot \nabla Z + \nabla X \cdot \nabla_Y \nabla Z
- \nabla_X \nabla Y \cdot \nabla Z - \nabla Y \cdot \nabla_X \nabla Z
-\nabla \nabla_Y X + \nabla \nabla_X Y \cdot \nabla Z.$$
Again, using Lemma \ref{pd}, $[\nabla_Y,\nabla]$ is order $1$, so that
$$\| [\nabla_Y,\nabla]X \cdot \nabla Z\|_{s-1} \le C \|X\|_s\|Y\|_s\|Z\|_s,$$
with the same estimate for $[\nabla_X, \nabla]Y \cdot \nabla Z$.  After 
commutation, most of the terms in (\ref{A}) cancel, and we are left to
estimate
$$ \hat{\triangle}^{-1} \nabla^* [ \nabla X \cdot \nabla_Y \nabla Z - 
\nabla Y \cdot \nabla_X\nabla Z].$$
It suffices to estimate the first term.  Now
\begin{equation}\label{a1}
\hat{\triangle}^{-1} \nabla^* [ \nabla X \cdot \nabla_Y \nabla Z]
= \hat{\triangle}^{-1} [ (\nabla_Y \nabla Z)^t \cdot \hat{\triangle}X^t]
+ \hat{\triangle}^{-1}(\nabla^* \nabla_Y \nabla Z),
\end{equation}
so the first term in the right-hand-side of (\ref{a1}) is clearly a 
continuous mapping in $H^s$.  For the second
term we use the identity on divergence-free vector fields given by
$$ \text{div} \nabla_X Y = \text{Ric}(X,Y) + \text{Tr}(\nabla X \cdot \nabla Y),$$
where $\text{Ric}(X,Y)= \langle Ric\langle X \rangle, Y\rangle$.  We obtain that
$$ \nabla^* \nabla_Y \nabla Z =
\text{grad}[ \text{Ric}(Y,Z) + \text{Tr}(\nabla Y \cdot \nabla Z) ] +
[\nabla^*,\nabla] \nabla_YZ + \nabla^* [\nabla_Y,\nabla]Z.$$
Hence, using Lemma \ref{pd}, $\nabla^* \nabla_Y \nabla Z:H^s \rightarrow
H^{s-2}$ is continuous, so that
$$ \| \hat{\triangle}^{-1} \nabla^* [ \nabla X \cdot \nabla_Y \nabla Z
- \nabla Y \cdot \nabla_X \nabla Z] \|_s \le C \|X\|_s \|Y\|_s \|Z\|_s.$$

This completes the estimates on each term  of $R^1_e(X,Y)Z$.  Since
we allow our constant to depend on $\eta$ and since $H^{s-2}$ is a
multiplicative algebra, we have that for any $\eta \in {\mathcal D}^s(M)$,
$$ \|R^1(X_\eta,Y_\eta)Z_\eta\|_s \le C \|X_\eta\|_s \|Y_\eta\|_s \|Z_\eta\|_s,$$
where $C$ denotes any constant which may depend on $s, \eta$, and derivatives
of $\langle \cdot, \cdot \rangle$ on $M$.
\end{proof}

\subsection{Curvature of $\tilde{\nabla}^1$}
Next, we define the (weak) curvature $\tilde{R}^1$ of the induced metric
$\langle \cdot, \cdot \rangle_1$ on ${\mathcal D}_\mu^s(M)$ as
\begin{equation}\nonumber
\begin{array}{c}
\tilde{R}^1_\eta: T_\eta{\mathcal D}_\mu^s(M) \times T_\eta{\mathcal D}_\mu^s(M)
\times T_\eta{\mathcal D}_\mu^s(M) \rightarrow T_\eta{\mathcal D}_\mu^s(M),\\
\tilde{R}^1_\eta(X_\eta,Y_\eta)Z_\eta = (\tilde{\nabla}^1_X \tilde{\nabla}^1_Y Z)_\eta -
(\tilde{\nabla}^1_Y \tilde{\nabla}^1_X Z)_\eta - (\tilde{\nabla}^1_{[X,Y]}Z)_\eta,
\end{array}
\end{equation}
where $\eta\in {\mathcal D}_\mu^s(M)$, and $X,Y,Z$ are smooth extensions of
$X_\eta,Y_\eta,Z_\eta$ in a neighborhood of $\eta$.

In order to estimate $\tilde{R}^1$, we shall make use of the Gauss formula
in submanifold geometry which relates the curvature of ${\mathcal D}^s(M)$
with the curvature of ${\mathcal D}_\mu^s(M)$ using the second fundamental
form.  Let $X,Y,Z$, and $W$ be smooth vector fields on ${\mathcal D}_\mu^s(M)$.
Then for any $\eta \in {\mathcal D}_\mu^s(M)$, we have
\begin{eqnarray}
\langle \tilde{R}^1(X,Y)Z,W \rangle_1 &=& \langle R^1(X,Y)Z,W \rangle_1 
+ \langle S_\eta(Y,Z), S_\eta(X,W) \rangle_1 \nonumber\\
&& \ \ \qquad \qquad \qquad \qquad - \langle S_\eta(X,Z), S_\eta(Y,W) \rangle_1 .
\label{gauss}
\end{eqnarray}

\begin{thm} \label{thmRu1}
The curvature $\tilde{R}^1$ of the induced $H^1$ metric on ${\mathcal D}_\mu^s(M)$
is a trilinear operator which is continuous in the $H^s$ topology for $s >
{\frac{n}{2}}+2$.
\end{thm}
\begin{proof}
For the purpose of obtaining estimates on $\tilde{R}^1$ we shall use the
equivalent $H^s$ metric given at the identity for $X,Y \in 
T_e {\mathcal D}_\mu^s(M)$ by
$$\langle X, Y \rangle_s = \langle X, (1-\hat{\triangle})^s Y\rangle_{L^2},$$
and then extended to $T{\mathcal D}_\mu^s(M)$ by right invariance.  This
gives a smooth invariant metric on ${\mathcal D}_\mu^s(M)$ which induces
a topology which is equivalent to the underlying topology of 
${\mathcal D}_\mu^s(M)$.

We will estimate $\sup_{\|W\|_s=1}\langle \tilde{R}^1(X,Y)Z,W \rangle_s$ 
using the Gauss formula (\ref{gauss}). Let 
$X,Y,Z \in T_e {\mathcal D}_\mu^s(M)$, and let $W\in C^\infty(TM)$, div$W=0$.
We have that
\begin{eqnarray}
\langle \tilde{R}^1(X,Y)Z,(1-\hat{\triangle})^sW \rangle_0 &= &
\langle R^1(X,Y)Z,(1-\hat{\triangle})^sW \rangle_0 \label{gauss2}\\
&& + \langle S_e(Y,Z), (1-\hat{\triangle})S_e(X,(1-\hat{\triangle})^{s-1}W)
 \rangle_0 \nonumber\\
&& \qquad - \langle S_e(X,Z), (1-\hat{\triangle})S_e(Y,
(1-\hat{\triangle})^{s-1}W) \rangle_0 . \nonumber
\end{eqnarray}

Now, $S_e(X,Y) = Q_e(\nabla_XY) + Q_eD(X,Y)$, where $D(X,Y)=A(X,Y)+
B(X,Y)+C(X,Y)$, so
\begin{eqnarray}
&&\langle S_e(Y,Z), (1-\hat{\triangle})S_e(X,(1-\hat{\triangle})^{s-1}W)
 \rangle_0\nonumber\\ 
&& \qquad \qquad =\langle Q_e(\nabla_YZ), 
(1-\hat{\triangle})Q_e\nabla_X(1-\hat{\triangle})^{s-1}W\rangle_0\nonumber\\
&& \qquad \qquad \qquad +
\langle Q_e(\nabla_YZ), (1-\hat{\triangle}) Q_e D(X,(1-\hat{\triangle})^{s-1}W)
\rangle_0\label{10} \\
&& \qquad \qquad \qquad \qquad +
\langle Q_eD(Y,Z), (1-\hat{\triangle}) Q_e (\nabla_X(1-\hat{\triangle})^{s-1}W)
\rangle_0 \nonumber\\
&& \qquad \qquad \qquad \qquad \qquad +
\langle Q_eD(Y,Z), (1-\hat{\triangle}) Q_e D(X,(1-\hat{\triangle})^{s-1}W)
\rangle_0 \nonumber.
\end{eqnarray}

For the first step, we will obtain the estimates for (\ref{10}) in the case
where $D$ is just $B$.
We begin by estimating the first term on the right-hand-side of (\ref{10}).
Using the fact that $Q_e$ is also an orthogonal projection in $L^2$, we have
that
\begin{eqnarray}
&&\langle Q_e(\nabla_YZ),(1-\hat{\triangle})Q_e\nabla_X(1-\hat{\triangle})
^{s-1}W\rangle_0 \label{11}\\
&&\qquad \qquad \qquad \qquad \qquad = 
-\langle (1-\hat{\triangle})^{\frac{s-2}{2}} \nabla_X
Q_e (1-\hat{\triangle}) Q_e \nabla_YZ, (1-\hat{\triangle})^{\frac{s}{2}}W
\rangle_0.\nonumber
\end{eqnarray}
Using the identity for divergence-free vector fields
$$\text{div} \nabla_XY = \text{Ric}(X,Y) + \text{Tr}(\nabla X \cdot \nabla Y),$$
and choosing a smooth local orthonormal frame $\{e_i\}$ in which the rough
Laplacian $\hat{\triangle}= \nabla_{e_i} \nabla_{e_i}$, we see that
\begin{eqnarray}
Q_e \hat{\triangle} Q_e \nabla_YZ &=& \text{grad} \hat{\triangle}^{-1}
\text{Ric}(e_i, \nabla_{e_i}\text{grad} \hat{\triangle}^{-1} \text{Ric} (Y,Z)) \nonumber\\
&+& \text{grad} \hat{\triangle}^{-1}\text{Tr}[\nabla e_i \cdot
\nabla \nabla_{e_i} \text{grad} \hat{\triangle}^{-1}\text{Ric}(Y,Z)]\nonumber\\
&+& \text{grad} \hat{\triangle}^{-1}
\text{Ric}(e_i, \nabla_{e_i}\text{grad} \hat{\triangle}^{-1} \text{Tr}
[\nabla Y \cdot \nabla Z] \nonumber\\
&+& \text{grad} \hat{\triangle}^{-1}\text{Tr}[\nabla e_i \cdot
\nabla \nabla_{e_i} \text{grad} \hat{\triangle}^{-1}\text{Tr}[\nabla Y
\cdot \nabla Z ]] \label{12}
\end{eqnarray}
We estimate the last term in (\ref{12}) since it is least regular.  We
obtain
\begin{eqnarray*}
&&\|\text{grad} \hat{\triangle}^{-1}\text{Tr}[\nabla e_i \cdot
\nabla \nabla_{e_i} \text{grad} \hat{\triangle}^{-1}\text{Tr}[\nabla Y
\cdot \nabla Z ]] \|_{s-1} \\
&&\qquad \qquad \qquad \le 
\|\text{Tr}[\nabla e_i \cdot
\nabla \nabla_{e_i} \text{grad} \hat{\triangle}^{-1}\text{Tr}[\nabla Y
\cdot \nabla Z ]] \|_{s-2} \\
&& \qquad \qquad \qquad \le C
\|\nabla_{e_i} \text{grad} \hat{\triangle}^{-1}\text{Tr}[\nabla Y
\cdot \nabla Z ]] \|_{s-1}
\end{eqnarray*}
where we used the fact that $H^{s-2}$ is a multiplicative algebra, and
the constant $C$ may depend on $e_i$.
Now
\begin{eqnarray*}
&&\|\nabla_{e_i} \text{grad} \hat{\triangle}^{-1}\text{Tr}[\nabla Y
\cdot \nabla Z ]] \|_{s-1}\\
&&\qquad \qquad \qquad \le 
\| \hat{\triangle}^{\frac{s-1}{2}}( \nabla \text{grad} \hat{\triangle}^{-1}
\text{Tr}[\nabla Y \cdot \nabla Z]) \cdot e_i \|_0 \\
&&\qquad \qquad \qquad \qquad +
\| (\nabla \text{grad} \hat{\triangle}^{-1}
\text{Tr}[\nabla Y \cdot \nabla Z]) \cdot \hat{\triangle}^{\frac{s-1}{2}}
e_i \|_0 \\
&&\qquad \qquad \qquad \le 
C \|\text{Tr}[\nabla Y \cdot \nabla Z] \|_{s-1} \le C \|Y\|_s \|Z\|_s.
\end{eqnarray*}
This shows that $\|Q_e (1-\hat{\triangle})Q_e \nabla_YZ\|_{s-1} \le C
\|Y\|_s \|Z\|_s$, so that applying the Cauchy-Schwartz inequality to
(\ref{11}) we obtain
\begin{eqnarray*}
&&|\langle Q_e(\nabla_YZ),(1-\hat{\triangle})Q_e\nabla_X(1-\hat{\triangle})
^{s-1}W\rangle_0|\\
&&\qquad \qquad \le 
C \|\nabla_X Q_e (1-\hat{\triangle}) Q_e \nabla_Y Z \|_{s-2} \|W\|_s\\
&&\qquad \qquad \le 
C\left\{\| \hat{\triangle}^{\frac{s-2}{2}}
( \nabla Q_e(1-\hat{\triangle})Q_e \nabla_YZ) \cdot X\|_0 \right. \\
&&\qquad \qquad \qquad \qquad \qquad \qquad \qquad \left.
+ \| \nabla Q_e(1-\hat{\triangle})Q_e \nabla_YZ \cdot \hat{\triangle}^
{\frac{s-2}{2}}X\|_0 \right\} \|W\|_s\\
&&\qquad \qquad \le 
C\left\{\| \nabla Q_e(1-\hat{\triangle})Q_e \nabla_YZ\|_{s-2} \|X\|_\infty \right. \\
&&\qquad \qquad \qquad \qquad \qquad \qquad \qquad \left.
+ \| \nabla Q_e(1-\hat{\triangle})Q_e \nabla_YZ\|_\infty 
\|X\|_{s-2} \right\} \|W\|_s\\
&&\qquad \qquad \le 
C \|Q_e (1-\hat{\triangle}) Q_e \nabla_Y Z \|_{s-1}\|X\|_s \|W\|_s\\
&&\qquad \qquad \le 
C \|X\|_s \|Y\|_s\|Z\|_s\|W\|_s.
\end{eqnarray*}
Since $B:H^s \times H^s \rightarrow H^{s+1}$ continuously, we have estimated
the first and third terms on the right-hand-side of (\ref{10}).

Next we estimate the second term on the right-hand-side of (\ref{10}).
We have that
\begin{eqnarray}
&&B(X,(1-\hat{\triangle})^{s-1}W) = {\frac{1}{2}}(1-\hat{\triangle})^{-1} 
\text{Tr} [R(\cdot,\nabla_\cdot (1-\hat{\triangle})^{s-1}W)X \label{13} \\
&& \qquad +[R(\cdot,\nabla_\cdot X)(1-\hat{\triangle})^{s-1}W
- \nabla_\cdot[R(X,
\cdot)(1-\hat{\triangle})^{s-1}W + R((1-\hat{\triangle})^{s-1}W,\cdot)X]]
\nonumber
\end{eqnarray}
Let us begin our estimate with the first of the four terms in
(\ref{13}).  Let 
$$V= {\frac{1}{2}}(1-\hat{\triangle})^{-1} Q_e (1-\hat{\triangle})
Q_e \nabla_YZ,$$
which is of Sobolev class $H^{s+1}$.  Then
\begin{eqnarray}
&&{\frac{1}{2}} |\langle Q_e(\nabla_YZ) , (1-\hat{\triangle}) Q_e 
(1-\hat{\triangle})^{-1} \text{Tr}
[R(e_i, \nabla_{e_i} (1-\hat{\triangle})^{s-1}W)X] \rangle_0|\nonumber \\
&& \qquad  =
|\langle V, \text{Tr}[R(e_i,\nabla_{e_i}(1-\hat{\triangle})^{s-1}W)X] 
\rangle_0.|\nonumber\\
&& \qquad  = \left|
\int_M \left\langle (1-\hat{\triangle})^{\frac{s-2}{2}}
\{ \langle \nabla_{e_i}V, R(e_i,X)\cdot\rangle^\sharp +
\langle V, (\nabla_{e_i}R)(e_i,X)\cdot\rangle^\sharp\right.\right.\nonumber\\
&& \qquad\qquad  +
\langle V, R(\nabla_{e_i}e_i,X)\cdot\rangle^\sharp +
\langle V, \text{div}(e_i) R(e_i,X)\cdot\rangle^\sharp \nonumber\\
&& \qquad\qquad \left.\left. +
\langle V, R(e_i,\nabla_{e_i}X)\cdot\rangle^\sharp \},
(1-\hat{\triangle})^{\frac{s}{2}} W \right\rangle dx\right|. \label{14}
\end{eqnarray}
Now
\begin{eqnarray*}
&&\left|\int_M \left\langle (1-\hat{\triangle})^{\frac{s-2}{2}}
\{ \langle \nabla_{e_i}V, R(e_i,X)\cdot\rangle^\sharp\},
(1-\hat{\triangle})^{\frac{s}{2}} W \right\rangle \right| \\
&&\qquad =\left|\int_M \left\langle 
\text{Tr}\left[ \langle(1-\hat{\triangle})^{\frac{s-2}{2}}
\nabla_\cdot V, R(\cdot,X)\rangle^\sharp +
\langle \nabla_\cdot V, (1-\hat{\triangle})^{\frac{s-2}{2}}
R(\cdot,X)\rangle^\sharp\right],\right. \right. \\
&& \qquad \qquad \qquad \qquad  \left. \left.
(1-\hat{\triangle})^{\frac{s}{2}} W \right\rangle \right| \\
&&\qquad \le 
\left\| \text{Tr} \left[ \langle
(1-\hat{\triangle})^{\frac{s-2}{2}} \nabla_\cdot V, R(\cdot,X)\rangle^\sharp
\right. \right. \\
&&\qquad \qquad  \left. \left.
+ \langle \nabla_\cdot V, ((1-\hat{\triangle})^{\frac{s-2}{2}}R)(\cdot,
X) + R(\cdot, (1-\hat{\triangle})^{\frac{s-2}{2}}X)\rangle^\sharp\right] \right\|_0
\|W\|_s\\
&&\qquad \le 
C\left[ \|\nabla V\|_{s-2} \|R\|_\infty \|X\|_\infty
+ \| \nabla V \|_\infty \| (1-\hat{\triangle})^{\frac{s-2}{2}}R\|_\infty\|X\|_\infty
\right. \\
&&\qquad \qquad  \left. 
+ \| \nabla V\|_\infty \|R\|_\infty \|X\|_{s-2}\right] \|W\|_s \\
&&\qquad \le C \|V\|_{s-1}\|X\|_s \|W\|_s \le C \|Q_e(1-\hat{\triangle})Q_e
\nabla_YZ\|_{s-3} \|X\|_s \|W\|_s\\
&&\qquad \le C \|X\|_{s}\|Y\|_s \|Z\|_s\|W\|_s,
\end{eqnarray*}
where the constant $C$ may depend on $M$,  the derivatives of the
metric $\langle \cdot, \cdot \rangle$ on $M$, and the local orthonormal
frame.  The remaining terms in (\ref{14}) can be estimated in the
same manner, so that
\begin{eqnarray*}
&&{\frac{1}{2}} |\langle Q_e(\nabla_YZ) , (1-\hat{\triangle}) Q_e 
(1-\hat{\triangle})^{-1} \text{Tr}
[R(e_i, \nabla_{e_i} (1-\hat{\triangle})^{s-1}W)X] \rangle_0|\nonumber \\
&& \qquad \le C \|X\|_{s}\|Y\|_s \|Z\|_s\|W\|_s.
\end{eqnarray*}

Using the same type of estimates, we may bound the remaining three terms
in (\ref{13}), so that the second term on the right-hand-side of (\ref{10})
with $D=B$
is majorized by $\|X\|_s\|Y\|_s\|Z\|_s\|W\|_s$.  The fourth term on
right-hand-side of (\ref{10}) with $D=B$ has more regularity than the second term,
and thus has the same majorization.  

Now, if we let $D=C$, we easily obtain the same estimates since $C$ is
as regularizing as $B$.  For $D=A$, we must estimate the term
$$ \langle Q_e \nabla_YZ, (1-\hat{\triangle}) Q_e (1-\hat{\triangle})^{-1}
\nabla^*( \nabla X \cdot \nabla (1-\hat{\triangle})^{\frac{s}{2}} W) \rangle_0.$$
With similar estimates as above, we can bound this term by
\begin{eqnarray*}
&&C \left( \|(1-\hat{\triangle})^{\frac{s-2}{2}} \text{grad div} X\|_0 \cdot
\| \nabla (1-\hat{\triangle})^{-1} Q_e (1-\hat{\triangle}) Q_e \nabla_Y Z\|_\infty
\right. \\
&& +  \| \text{grad div} X\|_\infty \cdot
\| (1-\hat{\triangle})^{\frac{s-2}{2}} \nabla (1-\hat{\triangle})^{-1} Q_e
(1-\hat{\triangle}) Q_e \nabla_Y Z \|_0 \\
&& + \| (1-\hat{\triangle})^{\frac{s-2}{2}} \nabla X\|_0 \cdot
\|\nabla ( \nabla (1-\hat{\triangle})^{-1}Q_e (1-\hat{\triangle}) Q_e
\nabla_YZ)^t \|_\infty \\
&& \left. + \|\nabla X\|_\infty \cdot
\|(1-\hat{\triangle})^{\frac{s-2}{2}} \nabla (1-\hat{\triangle})^{-1}Q_e 
(1-\hat{\triangle}) Q_e \nabla_YZ \|_0 \right)
\end{eqnarray*}
which is itself bounded by $C \|X\|_s\|Y\|_s\|Z\|_s\|W\|_s$.  The estimates for
the other terms involving $A$ are similar.

Hence, we have estimated the
second term on the right-hand-side of (\ref{gauss2}), and by symmetry
of the bound, the third term as well.  Proposition \ref{thmR1} gives us
the same majorization for the first term.

Since
\begin{eqnarray*}
\| \tilde{R}^1_e (X,Y)Z \|_s &=& \sup \{ \langle \tilde{R}^1_e
(X,Y)Z,W\rangle_s : \\
&& \qquad W \in C^\infty(TM), \text{div}W=0, \|W\|_s < 1 \} \\
&\le & C \|X\|_s\|Y\|_s\|Z\|_s,
\end{eqnarray*}
where $C$ depends on $M$ and the derivatives of the metric on $M$,
we have that $\tilde{R}^1_e$ is a bounded trilinear map on $H^s$.

Now the map $\eta \rightarrow P_\eta$ is continuously
differentiable, and since right translation only introduces terms of the
type $[T\eta]^{-1}$ and ${[T\eta]^{-1}}^t$, and as we have a multiplicative
algebra, the general case follows.
\end{proof}

\begin{rem}
One might try to argue that the boundedness in $H^s$ of 
$\tilde{R}^1$ follows immediately from the regularity of the geodesic
spray, but this argument fails for the following reason.  Let
${\mathcal U} \subset {\mathcal D}_\mu^s(M)$ be sufficiently small so
as to allow a trivialization of $T{\mathcal D}_\mu^s(M)$, and let
${\mathcal A}^1$ be the local connection $1$-form defining the 
$H^1$ covariant derivative $\tilde{\nabla}^1$.  The fact that
the geodesic spray of $\tilde{\nabla}^1$ is $C^1$ implies that
${\mathcal A}^1$ is a $C^1$ map as well.  Now the curvature can
be defined as $d {\mathcal A}^1 + {\mathcal A}^1 \wedge {\mathcal A}^1$,
and it may seem that for all $\eta \in {\mathcal U}$, 
$d {\mathcal A}^1(\eta)$ is then necessarily a continuous operator from
$H^s$ into $H^s$.  This is not the case, however, as the exterior derivative
is defined in terms of the $H^1$-Frechet derivative, while the fact
that ${\mathcal A}^1$ is $C^1$ is verified using the $H^s$-Frechet derivative.
It is for this reason, that curvatures of strong metrics are trivially bounded
operators in the strong topology of the manifold, while for weak metrics,
one must verify any boundedness claims.
\end{rem}

\subsection{Jacobi equations}
We can now prove the existence of solutions to the Jacobi equation
\begin{equation}\label{jacobi}
\tilde{\nabla}^1_{\dot{\eta}} \tilde{\nabla}^1_{\dot{\eta}} Y +
\tilde{R}^1_\eta(Y,{\dot{\eta}}){\dot{\eta}} =0
\end{equation}
along the geodesic $\eta(t)$ of the $H^1$-metric which solves the mean 
fluid motion
equation (\ref{spray}) in Lagrangian coordinates. Note that (\ref{spray})
may equivalently be written as
\begin{equation}\label{15}
\tilde{\nabla}^1_{\dot{\eta}} {\dot{\eta}}=0,
\end{equation}
for $\eta(t)$ a curve in ${\mathcal D}_\mu^s(M)$.  The Jacobi
equation (\ref{jacobi}) is the linearization  of (\ref{15}) along
the geodesic.

\begin{thm}
Let $s>{\frac{n}{2}}+2$ and let $Y_e,\dot{Y}_e \in T_e{\mathcal D}_\mu^s(M)$.
Then there exists a unique $H^s$ vector field $Y(t)$ along $\eta$ that
is a solution to (\ref{jacobi}) with initial conditions $Y(0)=Y_e$ and
$\tilde{\nabla}^1_{\dot{\eta}}Y(0)=\dot{Y}_e$.
\end{thm}
\begin{proof}
Let $\tau_t:T_e {\mathcal D}_\mu^s(M) \rightarrow 
T_{\eta(t)}{\mathcal D}_\mu^s(M)$ be the parallel translation along
$\eta$ induced by $\tilde{\nabla}^1$.  It is standard that
$\tau_t$ is a linear isomorphism such that $[\tau_t, \tilde{\nabla}^1]=0$,
and $\tau_t^*\langle\cdot,\cdot\rangle_1 = \langle \cdot, \cdot \rangle_1$.
We consider the curve in the algebra $V(t)=\tau_t^{-1}Y(t)$ where
$(d/dt)V(t)=\tau_t^{-1} \tilde{\nabla}^1_{\dot{\eta(t)}}Y(t)$, wherein
the Jacobi equation takes the form
$$\frac{d^2}{dt^2} V(t) = - \tau_t^{-1}\tilde{R}^1_{\eta(t)}(
\tau_t V(t), {\dot{\eta}}(t)) {\dot{\eta}}(t).$$
By Theorem \ref{thmRu1}, $\tilde{R}^1$ is bounded in $H^s$, so existence
and uniqueness immediately follow.
\end{proof}

\section{Stability and Curvature}
\label{sec5}
In this section, we define the notion of Lagrangian linear stability (see 
\cite{M1}).

\subsection{Lagrangian stability}
For $k \ge 1$, a fluid motion $\eta$ is Lagrangian $H^k$ (linearly)
stable if every solution of the Jacobi equation (\ref{jacobi}) along
$\eta$ is bounded in the $H^k$ norm.

\begin{thm}\label{thm_stab}
If $\eta(t)$ is a geodesic of $\tilde{\nabla}^1$ on ${\mathcal D}^s_\mu(M)$
whose pressure function $p(t)$ is constant for all $t$ and if the sectional 
curvature of $R^1$ is nonpositive, then $\eta$ is $H^k$ Lagrangian unstable
for $k \ge 1$.
\end{thm}
\begin{proof}
Let $\eta$ solve $\tilde{\nabla}^1_{\dot{\eta}}{\dot{\eta}} =0$ on 
${\mathcal D}_\mu^s(M)$, and let $Y(t)$ be a nontrivial Jacobi
field along $\eta$ with $Y(0)=0$, $\tilde{\nabla}^1_{\dot{\eta}} Y(0)=
\dot{Y}_e$.  If the sectional curvature of the plane spanned by
$Y(t)$ and ${\dot{\eta}}$ is nonpositive for $t$, then $\eta$ is
$H^k$ Lagrangian unstable for $k\ge 1$.  This follows from
Lemma 4.2 of \cite{M1} by replacing the $L^2$ norm with the $H^1$ norm.
Namely, for $t>0$, let $Z=Y/\|Y\|_1$ and compute
$$\tilde{\nabla}^1_{\dot{\eta}}\tilde{\nabla}^1_{\dot{\eta}}Y=
\frac{d^2}{dt^2}(\|Y\|_1)Z + 2\frac{d}{dt}(\|Y\|_1)
\tilde{\nabla}^1_{\dot{\eta}} Z +
\|Y\|_1
\tilde{\nabla}^1_{\dot{\eta}}\tilde{\nabla}^1_{\dot{\eta}}Y.$$
Taking the inner product of
$\tilde{\nabla}^1_{\dot{\eta}}\tilde{\nabla}^1_{\dot{\eta}}Y$ with
$Z$, and noting that $\|Z\|_1=1$ and that $Y$ solves (\ref{jacobi}),
we obtain that
$$\frac{d^2}{dt^2}(\|Y\|_1) = \left[
\| \tilde{\nabla}^1_{\dot{\eta}} Z \|_1^2 -
\langle \tilde{R}^1(Z,\dot{\eta})\dot{\eta},Z\rangle_1 \right]\|Y\|_1.$$
Thus, $(d^2/dt^2)\|Y\|_1 \ge 0$, so that $\|Y\|_1 > c t$ for all
$t >0$ and some positive constant $c$ depending on $\dot{Y}_e$, which
implies that $\|Y\|_k$ is unbounded for $k\ge 1$ by the compact
embedding: $H^k \hookrightarrow H^1$.

Since $\eta$ is a geodesic in ${\mathcal D}_\mu^s(M)$, Theorem
\ref{thm_spray} asserts that $U ={\dot{\eta}} \circ \eta^{-1}$ 
satisfies equation (\ref{7}) on $M$.  Thus, we have that
\begin{eqnarray*}
S_\eta({\dot{\eta}},{\dot{\eta}}) &=& 
       Q_\eta(\nabla^1_{\dot{\eta}} {\dot{\eta}})\\
&=&Q_e\left\{
\partial_t U + (1-\hat{\triangle})^{-1}\left[\nabla_{U}
(1-\hat{\triangle})U - \langle\nabla U \langle \cdot \rangle, 
\hat{\triangle} U \rangle^\sharp\right] 
\right\} \circ \eta\\
&=& -(\text{grad}\ p) \circ \eta = 0,
\end{eqnarray*}
so $\eta$ is a pressure constant geodesic of the right invariant $H^1$
metric on ${\mathcal D}_\mu^s(M)$ if and only if
$S_\eta({\dot{\eta}},{\dot{\eta}})=0$.

From the Gauss equation (\ref{gauss}),
$$
\langle \tilde{R}^1_\eta(X,{\dot{\eta}}){\dot{\eta}},X \rangle_1 =
\langle R^1(X,{\dot{\eta}}){\dot{\eta}},X)_1 - \|S_\eta({\dot{\eta}},X)\|^2_1,
$$
for any vector field $X(t)$ along the pressure constant geodesic $\eta$.  Hence,
$\langle \tilde{R}^1_\eta(X,{\dot{\eta}}){\dot{\eta}},X \rangle_1$
is nonpositive whenever
$\langle R^1(X,{\dot{\eta}}){\dot{\eta}},X\rangle_1$ is nonpositive.
\end{proof}

\begin{rem}
Note that on the flat torus ${\mathbb T}^n$,
the formula (\ref{3}) simplifies to $\nabla^1_XY = \nabla^0_XY
+A(X,Y)$, and since $R^0 = 0$, we have that for $X,Y,Z \in T_e{\mathcal D}_\mu^s(M)$,
\begin{eqnarray}
R^1_e(X,Y)Z &=& A_e(X, \nabla^1_YZ) - A_e(Y, \nabla^1_XZ) 
+\nabla^0_X A_e(Y,Z) - \nabla^0_Y A_e(X,Z) \nonumber\\
 && + A_e(X,A_e(Y,Z))- A_e(Y,A_e(X,Z))- A_e([X,Y],Z). \label{16}
\end{eqnarray}
Choose a coordinate chart $(U,x^i)$ on $M$.  At the identity $e$, 
$$2 A_e(X,Z)= (1-\triangle)^{-1}[\nabla^*( \nabla X \cdot \nabla Z
+ \nabla Z \cdot \nabla Z)].$$
Substitution of
$(1-\triangle)^{-1} \nabla^* (\nabla X \cdot \nabla Z)$ into
(\ref{16}) yields
\begin{eqnarray*}
&& \frac{\partial}{\partial x^j} \left[ (1-\triangle)^{-1}
\frac{\partial}{\partial x^l} \left( \frac{\partial Y^l}{\partial x^i}
\frac{\partial Z^i}{\partial x^n}\right) \right] X^j -
 \frac{\partial}{\partial x^l} \left[ (1-\triangle)^{-1}
\frac{\partial}{\partial x^j} \left( \frac{\partial X^j}{\partial x^i}
\frac{\partial Z^i}{\partial x^n}\right) \right] Y^l \\
&&+
(1-\triangle)^{-1} \frac{\partial }{\partial x^j}
\left[ \frac{\partial X^j}{\partial x^i}\frac{\partial }{\partial x^n}
\left(\frac{\partial Z^i}{\partial x^l} Y^l \right)\right]
-
(1-\triangle)^{-1} \frac{\partial }{\partial x^l}
\left[ \frac{\partial Y^l}{\partial x^i}\frac{\partial }{\partial x^n}
\left(\frac{\partial Z^i}{\partial x^j} X^j \right)\right]\\
&&+
(1-\triangle)^{-1} \frac{\partial }{\partial x^j}
\left\{ \frac{\partial Y^j}{\partial x^n}\frac{\partial }{\partial x^k}
\left[(1-\triangle)^{-1}
\frac{\partial }{\partial x^l}
\left(\frac{\partial X^l}{\partial x^i}\frac{\partial Z^i}{\partial x^n}
\right)\right]\right\} \\
&&\qquad \qquad -
(1-\triangle)^{-1} \frac{\partial }{\partial x^l}
\left\{ \frac{\partial Y^l}{\partial x^n}\frac{\partial }{\partial x^k}
\left[(1-\triangle)^{-1}
\frac{\partial }{\partial x^j}
\left(\frac{\partial X^j}{\partial x^i}\frac{\partial Z^i}{\partial x^n}
\right)\right]\right\}.
\end{eqnarray*}

It is clear that $R^1_e$ vanishes when $X,Y,Z$ have components of the form 
$e^{i\langle k,x\rangle}$.  More interestingly, one may compute the
sectional curvature $\langle R^1_e(X,Y)Y,X \rangle_1$ in the directions
$X=\sin (\langle k, x\rangle)\frac{\partial }{\partial x^1}+
\cos (\langle m, x\rangle)\frac{\partial }{\partial x^2}$ and
$Y=\cos (\langle k, x\rangle)\frac{\partial }{\partial x^1}+
\sin (\langle m, x\rangle)\frac{\partial }{\partial x^2}$.
For example, when $X=(\sin(kx^1),0)$ and $Y=(0,\cos(kx^2)$,
$$\langle R^1_e(X,Y)Y,X \rangle_1=0,$$ 
whereas if $X=(\sin(kx^1),0)$ and $Y=(\cos(kx^1),0)$, then
$$\langle R^1_e(X,Y)Y,X \rangle_1 < 0$$
for any choice of $k\neq 0$ (cf. \cite{M3}).
Recall that this computation of the curvature tensor of the
full diffeomorphism group is restricted to divergence free
vector fields, since we are ultimately only interested in the
stability of the motion on the volume preserving subgroup.  
\end{rem}

If $\eta$ is a geodesic in ${\mathcal D}_\mu^s(M)$, two points $\eta(t_1)$
and $\eta(t_2)$ are conjugate with respect to $\eta$ if there exists a nonzero
Jacobi field $Y(t)$ along $\eta$ such that $Y(t_1)=Y(t_2) = 0$.  Such Jacobi
fields are thus stable perturbations of the initial flow.  

\begin{cor}
Let $\eta$ be a pressure constant geodesic in ${\mathcal D}_\mu^s(M)$. If
the sectional curvature of $R^1$ is nonpositive, then there are no
conjugate points along $\eta$.
\end{cor}

\subsection{Examples}
\begin{exam}
A trivial example of a pressure constant geodesic in 
${\mathcal D}_\mu({\mathbb T}^2)$ is given by
$$ \eta(t)(x^1,x^2) = (x^1+h(x^2), x^2+ct),$$
where $c$ is a constant and $h$ is a smooth periodic function.  Let
\begin{eqnarray*}
G(\eta)&=& -D(\dot\eta \circ \eta^{-1})^tD(\dot\eta \circ \eta^{-1})
{[T\eta]^{-1}}^t + D(\dot\eta \circ \eta^{-1})D(\dot\eta \circ \eta^{-1})
{[T\eta]^{-1}}^t\\
&&+D(\dot\eta \circ \eta^{-1})D(\dot\eta \circ \eta^{-1})^t
{[T\eta]^{-1}}^t.
\end{eqnarray*}
Then on ${\mathbb T}^n$, equation (\ref{15}) simplifies to 
$$\ddot{\eta}\circ\eta^{-1}-\text{grad} \triangle^{-1} \text{Tr}[D({\dot{\eta}} \circ
\eta^{-1})]^2 = (\text{Id}-\text{grad}\triangle^{-1}\text{div})[
(1-\hat{\triangle}_\eta)^{-1}G(\eta)],$$
and since ${\dot{\eta}}(x^1,x^2)=(0,c)$, then $\eta$ is a geodesic.
\end{exam}

\begin{exam}
Another example of a pressure constant geodesic  in
${\mathcal D}_\mu({\mathbb T}^2)$ is given by
$$ \eta(t)(x^1,x^2) = (x^1+th(x^2), x^2),$$
where again $c$ is a constant and $h$ is a smooth periodic function.  In
this case
$$ {\dot{\eta}}\circ \eta^{-1}(y^1,y^2)=(h(y^2),0),$$
and we must verify that
\begin{eqnarray}
0&=& P_e \circ \left\{ \partial_t(\dot{\eta}\circ \eta^{-1}) + 
(1-\hat{\triangle})^{-1}\left[
\nabla_{{\dot{\eta}}\circ \eta^{-1}} (1-\hat{\triangle})({\dot{\eta}} \circ
 \eta^{-1}) \right.\right. \nonumber\\
&&\left.\left. - [\nabla {\dot{\eta}} \circ \eta^{-1}]^t \cdot 
\hat{\triangle}({\dot{\eta}} \circ \eta^{-1})\right] \right\}.
\label{e2}
\end{eqnarray}

Notice that for our choice of $\eta$, 
$(1-\hat{\triangle})^{-1} [\nabla U]^t \cdot \hat{\triangle} U =
\text{grad} F$, for some $F \in C^\infty(M)$; hence,
$P_e \circ (1-\hat{\triangle})^{-1} [\nabla U]^t \cdot \hat{\triangle} U =0$,
so that (\ref{e2}) is simply
\begin{equation}\label{e3}
\partial_t({\dot{\eta}}\circ \eta^{-1}) +
(1-\hat{\triangle})^{-1}
\nabla_{{\dot{\eta}}\circ \eta^{-1}} (1-\hat{\triangle})({\dot{\eta}} \circ
 \eta^{-1}) = -\mathrm{grad }\  p.
\end{equation}
But the left-hand-side of (\ref{e3}) vanishes, so $\eta$ is a pressure
constant geodesic.
\end{exam}
\begin{rem}
Theorem \ref{thm_stab} and the remarks which follow its proof imply that
the geodesic flows of the previous two examples with $h(x^2)= 
\sin (kx^2)$ are unstable to perturbations in the $\cos (k x^2)$ direction.
Other such examples of unstable perturbations can be constructed.
\end{rem}

\section*{Acknowledgments}
The author would like to gratefully acknowledge  and thank Darryl Holm, 
Jerry Marsden, Gerard Misio\l ek, and Tudor Ratiu for their continued interest and
support, as well as the many stimulating discussions and helpful suggestions they have
provided.

\end{document}